\documentclass[11pt,reqno]{amsart}
\usepackage{graphicx}
\usepackage{amsfonts, amsmath, amstext, amssymb, amscd, color, rotating}
\usepackage{hyperref} 
\usepackage{enumerate}
\usepackage[ textheight=615pt, textwidth=360pt]{geometry}
\newtheorem{theorem}{Theorem}[section] 
\newtheorem{lemma}[theorem]{Lemma}    
\newtheorem{corollary}[theorem]{Corollary}
\newtheorem{proposition}[theorem]{Proposition}

\title[Carath\'{e}odory balls and proper holomorphic maps]{Carath\'{e}odory balls and proper holomorphic maps on multiply-connected planar domains}

\author{Tuen Wai Ng}
\address{The University of Hong Kong, Pokfulam, Hong Kong}
\email{ntw@maths.hku.hk}

\author{Chiu Chak Tang }
\address{The University of Hong Kong, Pokfulam, Hong Kong}
\email{ChiuChakTang@gmail.com}

\author{Jonathan Tsai}
\address{The University of Hong Kong, Pokfulam, Hong Kong}
\email{jonathan.tsai@cantab.net}

\allowdisplaybreaks
\begin{document}
\maketitle
\begin{abstract}
In this paper, we will establish the inequivalence of closed balls and the closure of open balls under the Carath\'{e}odory metric in some planar domains of finite connectivity greater than $2$, and hence resolve a problem posed by Jarnicki, Pflug and Vigu{\'e} in 1992. We also  establish a corresponding result for some pseudoconvex domains in $\mathbb{C}^n$ for $n \ge 2$.

This result will follow from an explicit characterization (up to biholomorphisms) of proper holomorphic maps from a non-degenerate finitely-connected planar domain, $\Omega$, onto the standard unit disk $\mathbb{D}$ which answers a question posed by Schmieder in 2005. Similar to Bell and Kaleem's characterization of proper holomorphic maps in terms of Grunsky maps (2008), our characterization of proper holomorphic maps from $\Omega$ onto $\mathbb{D}$ is an analogous result to Fatou's famous result that proper holomorphic maps of the unit disk onto itself are finite Blashcke products. 

Our approach uses a harmonic measure condition of Wang and Yin (2017) on the existence of a proper holomorphic map with prescribed zeros. We will see that certain functions $\eta(\cdot,p)$ play an analogous role to the M\"{o}bius transformations that fix $\mathbb{D}$ in finite Blaschke products. These functions $\eta(\cdot,p)$ map $\Omega$ conformally onto the unit disk with circular arcs (centered at 0) removed and map $p$ to $0$ and can be given in terms of the Schottky-Klein prime function. We also extend a result of Grunsky (1941) and hence introduce a parameter space for proper holomorphic maps from $\Omega$ onto $\mathbb{D}$.
\end{abstract}
\section{Introduction} 
\label{intro}
For a domain $\Omega\subsetneq \mathbb{C}$ and any $z,z_{0}\in\Omega$, following \cite{jarnickipflug}, one can define the {\it Carath\'{e}odory pseudodistance} between $z$ and $z_{0}$ to be $c_{\Omega}(z,z_{0})=\tanh^{-1}(c^{*}_{\Omega}(z,z_{0}))$, where
\[c^{*}_{\Omega}(z,z_{0})=\sup\{|f(z)|: f\in\mathcal{H}(\Omega),f(z_{0})=0\}\]
is the {\it M\"obius pseudodistance} and
$\mathcal{H}(\Omega)$ denotes the set of holomorphic maps of $\Omega$ into the unit disc $\mathbb{D}$. The 
pseudodistance $c_{\Omega}$ was introduced by  
Carath\'{e}odory \cite{Caratheodory27} in 1927. It is known that $c_{\Omega}$ (as well as $c^{*}_{\Omega}$) is a distance if $\Omega$ is biholomorphic to a bounded domain in $\mathbb{C}$ and in this case, the topology induced by $c_{\Omega}$ and $c^{*}_{\Omega}$ are the same as that of the Euclidean topology of $\Omega$ (see Chapter 2 of \cite{jarnickipflug}). In general, there is a significant difference between the topological and geometrical structure of the space $(\Omega, c_{\Omega})$ (see \cite{jarnicki1991caratheodory}) and even basic 
properties of its open balls are still unknown. Let $B_{c_{\Omega}}(a;r)=\{z\in\Omega: c_{\Omega}(z,a)<r\}$ and $\overline{B}_{c_{\Omega}}(a;r)=\{z\in\Omega:c_{\Omega}(z,a)\le r\}$ denote the open ball and closed ball with the center at $a$ and radius $r>0$ in the Carath\'{e}odory metric respectively. It is known that $c_{\Omega}$ is continuous and hence $\overline{B_{c_{\Omega}}(a;r)} \subset \overline{B}_{c_{\Omega}}(a;r)$ (the closure is taken in the sense of the Euclidean topology of $\Omega$). It is then natural to ask if the closure of each open ball is the closed ball of the same centre and radius. This problem was first posed by 
Jarnicki, Pflug and Vigu{\'e} \cite{jarnicki1992} in 1992 and then repeated in Jarnicki and Pflug's book \cite{jarnickipflug} (see p.42 in the first edition and p. 66 in the second edition). As we will see later, it is also related to the problem of whether the Carath\'{e}odory balls for planar domains are connected or not. We will solve these two problems by proving the following result.
\begin{theorem}\label{thm:cball}
For any $g\geq 2$, there exists a  domain $\Omega \subset \mathbb{C}$ of connectivity $g+1$ such that for some $0<r_2<r_1$ and $z_0\in \Omega$, we have
\begin{enumerate}
\item[i)]
$B_{c_{\Omega}}(z_{0};r_{1})$ is disconnected and relatively compact in $\Omega$; 
\item[ii)]  $\overline{B_{c_{\Omega}}(z_{0};r_{2})} \subsetneq \overline{B}_{c_{\Omega}}(z_{0};r_{2})$. 
\end{enumerate}
\end{theorem}
Notice that when $\Omega\subsetneq \mathbb{C}$ is simply connected, all the Carath\'{e}odory balls are connected as in this case, the Carath\'{e}odory metric is the same as the Poincar\'e metric. The case where $\Omega$ is doubly-connected was resolved by Schmieder \cite{schmieder} and Frerick and Schmieder \cite{schmieder2} and in this case, all Carath\'{e}odory balls are connected. Therefore, with Theorem \ref{thm:cball}, the problem posed by Jarnicki, Pflug and Vigu{\'e} \cite{jarnicki1992} is completely resolved.
\ \\

Theorem \ref{thm:cball} also allows us to construct domains in $\mathbb{C}^{n}$ ($n\ge 2$) which satisfy (i) and (ii) in the statement of the theorem. For any domain $G \subset \mathbb{C}^{n-1}$ and $a,z \in G, \lambda \in \mathbb{D}$, we have the product formula $c_{G \times\mathbb{D}}((a,0),(z,\lambda))=\max\{c_{G}(a,z),c_{\mathbb{D}}(0,\lambda)\}$ for the  Carath\'{e}odory pseudodistance of the product domain $G \times\mathbb{D}$. This implies that $B_{c_{G\times\mathbb{D}}}((a,0);r)=B_{c_{G}}(a;r)\times B_{c_\mathbb{D}}(0;r)$ and hence by Theorem \ref{thm:cball}, we have the following   

\begin{corollary}\label{cor:cball}
Let $n \ge 2$ and $\Omega$ be the domain in Theorem \ref{thm:cball}. Then $M=\Omega\times\mathbb{D}^{n-1}$ is a bounded pseudoconvex domain in $\mathbb{C}^n$ such that for some $0<r_{2}<r_{1}$ and $z_0\in M$, we have 
\begin{enumerate}
\item[i)]
$B_{c_{M}}(z_{0};r_{1})$ is disconnected and relatively compact in $M$; 
\item[ii)]  $\overline{B_{c_{M}}(z_{0};r_{2})} \subsetneq \overline{B}_{c_{M}}(z_{0};r_{2})$. 
\end{enumerate}
\end{corollary}

The above bounded domain $\Omega\times\mathbb{D}^{n-1}$ is pseudoconvex but not strongly pseudoconvex in $\mathbb{C}^n$. This is because any open subset of the complex plane is a domain of holomorphy and a product of two domains of holomorphy is again a domain of holomorphy. Moreover, a domain in $\mathbb{C}^n$ is pseudoconvex if and only if it is a domain of holomorphy. Finally, it is not difficult to show that $\Omega\times\mathbb{D}^{n-1}$ is not strongly pseudoconvex (cf. the arguments used for the poly-disks in Example 76 of \cite{slatyer2016levi}). So the bounded domains $\Omega\times\mathbb{D}^{n-1}$
in Corollary \ref{cor:cball} are new examples of complex manifolds which  
satisfy conditions (i) and (ii) above because so far the only known examples are the bounded strongly pseudoconvex domain in $\mathbb{C}^{n}$ ($n \ge 2$) and  the non-planar hyperbolic Riemann surface constructed by Jarnicki, Pflug and Vigu{\'e} in \cite{jarnicki1992} (see also Chapter 2 of \cite{jarnickipflug}). 
\ \\

Returning to the case when $\Omega\subsetneq \mathbb{C}$,
Grunsky showed (in \cite{grunsky1,grunsky2}) that the supremum in the definition of $c^{*}_{\Omega}(z,z_{0})$ is attained by a proper holomorphic mapping of $\Omega$ onto $\mathbb{D}$ whose degree is equal to the connectivity of $\Omega$.
Recall that, a function $f:X \to Y$ from a topological space $X$ to a topological space $Y$ is said to be \textit{proper} if $f^{-1}(K)$ is compact for any compact set $K \subset Y$. Theorem \ref{thm:cball} will follow from an explicit characterization of proper holomorphic maps mapping a non-degenerate finitely-connected domain $\Omega$ onto $\mathbb{D}$ which we establish in Theorems \ref{thm:expression} and \ref{thm:expression2} of this paper.  This characterization also answers Question 1 posed by Schmieder in \cite{schmieder}. 

Let $\mathbb{D}=\{ z \in \mathbb{C} \: : \: |z| <1 \}$ be the open unit disk in $\mathbb{C}$. In 1923, Fatou \cite{fatou} proved that $f:\mathbb{D} \to \mathbb{D}$ is a proper holomorphic map from $\mathbb{D}$ onto $\mathbb{D}$ if and only if $f$ is a finite Blaschke product, i.e., $f$ takes the form
\[ 
f(z) = e^{i \theta} \prod^n_{k=1} \dfrac{z-p_k}{1-\overline{p_k}z}
\]
where $\theta \in [0,2\pi)$ and $p_k \in \mathbb{D}$ for all $k=1,\ldots ,n$. 

Let $\mathbb{A}_r=\{ z \in \mathbb{C} \: : \: r<|z| <1 \}$ be an annulus in $\mathbb{C}$ and let $C_{r}=\{z\in\mathbb{C}:|z|=r\}$. In 2017, Wang and Yin \cite{Wang} proved that there exists a proper holomorphic map $f:\mathbb{A}_r \to \mathbb{D}$ if and only if the zeros of $f$, $p_{1},\ldots,p_{n}$ satisfy $|p_1 p_2 \cdots p_n|=r^d$ where  $d=\deg (\left. f \right|_{C_{r}} )$. Furthermore, the proper holomorphic map is given by
\[ 
f(z) = \dfrac{e^{i \theta}}{z^d} \prod^n_{k=1}
\left(
\dfrac{z-p_k}{1-\overline{p_k}z}
\prod^{\infty}_{j=1}
\dfrac{(z-p_kr^{2j})(z-p_kr^{-2j})}{(1-\overline{p_k}r^{2j}z)(1-\overline{p_k}r^{-2j}z)}
\right)
\]
for some $\theta \in [0,2\pi)$ 

From now on, let $g \ge 1$ and $\Omega$ be a non-degenerate $(g+1)-$connected domain in $\mathbb{C}$. Let $\gamma_0, \ldots ,\gamma_g$ be the boundaries of the connected components of the complement of $\Omega$ in the extended complex plane. The assumption that $\Omega$ is non-degenerate means that each boundary component $\gamma_{0},\ldots,\gamma_{g}$ is not a point.

 For each $l=0,\ldots , g$, let $u_l : \Omega \to \mathbb{R}$ denote the harmonic function on $\Omega$ satisfying $\left. u_l \right|_{\gamma_m} = \delta_{lm}$, where $\delta_{lm}$ is the Kronecker delta. In other words, $u_{l}(z)$ is the harmonic measure of $\gamma_{l}$ with respect to $z$ in $\Omega$. In \cite{Wang},  Wang and Yin proved that a proper holomorphic map $f: \Omega \to \mathbb{D}$ of degree $n$  with zeros $p_{1},\ldots,p_{n}$ exists if and only if for some $n_{0},n_{1},\dots, n_{g}\in\mathbb{N}$ with $n=n_{0}+\cdots+n_{g}$, 
\begin{equation}
\label{eq:condition1}
     \sum_{k=1}^n u_j (p_k) = n_j
     \qquad
     \mbox{for all $j=1,\dots , g$}
\end{equation}
and the degree of $f$ restricted to  $\gamma_l$ is $n_l$ for $l=0,\ldots,g$. Moreover, the map $f$ is unique up to rotation. Note also that in 2019, Bogatyrev \cite{bogatyrev2019} proved a necessary and sufficient condition for proper homomorphic maps from a genus $g$ Riemann surface with $k>0$ boundary components onto the closed unit disk $\overline{\mathbb{D}}$. He also gave a representation for such proper holomorphic maps in terms of the third-kind abelian differentials. 

In this paper, we will make use of Wang and Yin's condition (\ref{eq:condition1}) to give an explicit description of $\mathcal{P}(\Omega)$, the set of all proper holomorphic maps of $\Omega$ onto $\mathbb{D}$ when $\Omega$ is a {\it circular domain}, i.e, $\Omega$ is $\mathbb{D}$ excluding $g$ mutually disjoint closed disks inside $\mathbb{D}$. Notice that there is no loss of generality to assume that $\Omega$ is a  circular domain because by Koebe's generalised Riemann mapping theorem \cite{koebe1,koebe2}, any
non-degenerate finitely connected domain in $\mathbb{C}$ is biholomorphic to a circular domain. So from now on, $\gamma_{0}$ is the unit circle and $\gamma_{1},\ldots,\gamma_{g}$ are disjoint circles inside $\mathbb{D}$. In addition, a \textit{circularly slit disk }is the unit disk with circular arcs (centered at 0) removed. It is a classical result that there exists a conformal map from $\Omega$ onto a circularly slit disk which preserves the unit circle and maps $p\in\Omega$ to $0$ (see e.g. Chapter 15 in \cite{conway}). Moreover, this map is unique up to rotation. The following result gives an explicit description of $\mathcal{P}(\Omega)$,

\begin{theorem}
\label{thm:expression}
Suppose that $p_1, \ldots, p_n$ are points in $\Omega$ satisfying the condition (\ref{eq:condition1}) for positive integers $n_1,\ldots,n_g$. Then up to a multiplication factor $e^{i\theta}$, any proper holomorphic map $f\in\mathcal{P}(\Omega)$ which has zeros $p_1, \ldots, p_n$ and satisfies $\deg ( \left. f \right|_{\gamma_j}  )= n_j$ for $j=1,\ldots,g$ can be written uniquely as
\[
f(z)= \exp \left( -2 \pi i \sum\limits_{j=1}^{g} n_j v_j (z) \right) \prod_{k=1}^n \eta (z,p_k).
\]
Here $v_j$ denotes the integral of the first kind for $j=1,\ldots, g$ and $\eta( \cdot ,p)$ is the conformal map from $\Omega$ onto a circularly slit disk which maps $p$ to zero and preserving the unit circle. 
\end{theorem}
We will also give the following alternate formula for a proper holomorphic map along with an alternate condition on the zeros of the proper holomorphic map that guarantees existence of the proper holomorphic map with those zeros.
\begin{theorem}
\label{thm:expression2}
Suppose $P=\left\lbrace p_{l,k} : l=0,\ldots,g
\: \mbox{and} \:
 k=0,\ldots,n_{l}-1\right\rbrace
\subset\Omega$ is a set of $n$ points in $\Omega$ that satisfy the condition (\ref{eq:condition1}). Then the unique proper holomorphic maps whose zeros are $p_{l,k}$ for $l=0,\ldots,g$ and  $k=0,\ldots,n_{j}-1$ and with $f(1)=1$ can be written as
\[
f(z)= K(P) \prod_{l=0}^{g}\prod_{k=0}^{n_{l}-1} \eta_{l} (z,p_{l,k})
\]
where $K(P)$ is a constant that does not depend on $z$; 
$\eta_{0}(\cdot,p)$ maps $\Omega$ conformally onto a circularly slit disk such that $\partial\mathbb{D}$ is mapped to $\partial\mathbb{D}$ and $\eta_{0}(p,p)=0$; and for $j=1,\ldots,g$,
$\eta_{j}(\cdot,p)$ maps $\Omega$ conformally onto a circularly slit disk such that $\gamma_{j}$ is mapped to $\partial\mathbb{D}$ and $\eta_{j}(p,p)=0$. 

Moreover, condition (\ref{eq:condition1}) on the zeros  can be rewritten as:
\ \\
A proper holomorphic map with zeros $p_{1},\ldots,p_{n}$ exists if and only if we can reindex the zeros as 
$\left\lbrace p_{l,k} : l=0,\ldots,g
\: \mbox{and} \:
 k=0,\ldots,n_{l}-1\right\rbrace$
such that
\begin{equation}  \prod_{l=0}^{g}\prod_{k=0}^{n_{l}-1}\rho_{l,i}(p_{l,k}) \text{ does not depend on } i, \label{eq:condition3}\end{equation}
where $\rho_{l,i}(p_{l,k})$ denotes the radius of the image of $\gamma_{i}$ under $\eta_{l}(\cdot,p_{l,k})$.
\end{theorem}
Note that for $l=0,\ldots,g$ and $p\in\Omega$, the functions $\eta_{l}(\cdot,p)$ have constant modulus on each boundary component. So if a function $F$ is a product of functions of the form $\eta_{l}(\cdot,p)$ (for some $l=0,\ldots,g$ and $p\in\Omega$), then $F$ will also have constant modulus $R_{l}$ on each boundary component $\gamma_{l}$ for $l=0,\ldots,g$. The condition (\ref{eq:condition3}) means that if $R_{0}=\cdots=R_{g}$ (i.e. $F$ has constant modulus on $\partial\Omega$), then after scaling, $F$ is a proper holomorphic map.
\ \\

Theorems \ref{thm:expression} and \ref{thm:expression2} show that the conformal maps $\eta(\cdot,p)$ and $\eta_{j}(\cdot,p)$ play an analogous role to the Blaschke factors
\[m(z,p)\equiv\frac{z-p}{1-\overline{p}z}\]
in a finite Blaschke product. In \cite{crowdy13}, Crowdy and Marshall give an expression for the function $\eta(\cdot,p)$ and $\eta_{j}(\cdot,p)$ in terms of the Schottky-Klein prime function, $\omega(z,\zeta)$, of the domain $\Omega$ (note that $\eta = \eta_0$). In particular,
\begin{equation*}\eta(z,p)\equiv\frac{\omega(z,p)}{|p|\omega(z,\overline{p}^{-1})}.\end{equation*}
The Schottky-Klein prime function is constructed using the Schottky group $\Theta$, a discrete group of M\"{o}bius transformations that is associated to the domain $\Omega$. We will see later that for certain Schottky groups, we can write
\begin{equation}\eta(z,p)\equiv \prod_{\theta\in\Theta}\dfrac{m(\theta(z),p)}{m(\theta(1),p)}\label{eqn:embedding}.\end{equation}

 Let $\nu=(n_{0},n_{1},\ldots, n_{g}) \in \mathbb{N}^{g+1}$. We will call $\nu$, the \textit{boundary degree} of a proper holomorphic map $f \in \mathcal{P}(\Omega)$ if $\deg (\left. f \right |_{\gamma_l})= n_l$ for $l=0,\ldots,g$. Let $\mathcal{P}_{\nu}(\Omega)$ denote the family of proper holomorphic maps of boundary degree $\nu$ and $\mathcal{P}_{n}(\Omega)$ denote the family of proper holomorphic maps of degree $n$. In the special case where $\nu=(1,\ldots,1)$, it is clear that $\mathcal{P}_{\nu}(\Omega)=\mathcal{P}_{g+1}(\Omega)$. Note that the Riemann-Hurwitz formula implies that $g+1$ is the minimum possible degree of a proper holomorphic map of $\Omega$ onto $\mathbb{D}$. 
 
 The requirement that the zeros satisfy (\ref{eq:condition1}) means that we cannot easily construct a proper holomorphic map by prescribing points in $\Omega$ that are its zeros. A theorem of Grunsky \cite{grunsky1941ueber} states that we can uniquely specify a proper holomorphic map in $\mathcal{P}_{g+1}(\Omega)$ by prescribing one point in $\Omega$ to be its zero and one point on each of the $g+1$ boundaries of $\Omega$ to be the preimages of $1$ (see Section 6.5 in \cite{book_goluzin1969}).  For completeness, we state this theorem here, following Schmieder's formulation \cite{schmieder}.
 
 \begin{theorem}[Grunsky \cite{grunsky1941ueber}]
 Suppose that $p\in\Omega$ and
 $w_{l}\in\gamma_{l}$ for each $l=0,\ldots, g$.
There exists a unique proper holomorphic map  $f\in\mathcal{P}_{g+1}(\Omega)$ such that $f(p)=0$ and $f(w_{l})=1$ for each $l=0,\ldots, g$. 
 \end{theorem}
 
We will prove the following theorem which extends this result to $\mathcal{P}_{\nu}(\Omega)$ for any boundary degree $\nu$. We remark that this also provides a new proof of the above theorem.
\begin{theorem}
\label{thm:parameter}
Suppose that $p\in\Omega$ and $w_{l,k}\in\gamma_{l}$ for each $l=0,\ldots,g$, $k=0,\ldots, n_{l}-1$. 
Then there exists a proper holomorphic map  $f\in\mathcal{P}_{\nu}(\Omega)$ such that
\begin{itemize}
\item[(i)]  $f(p)=0$;
\item[(ii)] $f(w_{l,k})=1$ for each $l=0,\ldots,g$ and $k=0,\ldots, n_{l}-1$.
\end{itemize}
If $n=g+1$, such $f$ is unique. If $n>g+1$, then for each $l=0,\ldots,g$ such that $n_{l}> 1$ and $k=1,\ldots,n_{l}-1$, let $\lambda_{l,k}$ be a given positive number. Then there exists a unique proper holomorphic map  $f\in\mathcal{P}_{\nu}(\Omega)$ such that (i) and (ii) above are satisfied as well as the condition (iii) that
for each such $l$ and $k$, \begin{equation} \frac{|f'(w_{0,0})|}{|f'(w_{l,k})|}=\lambda_{l,k}.\label{eqn:param}\end{equation}
\end{theorem}

Notice that for any proper map $f:\Omega \to \mathbb{D}$, $f$ maps each $\gamma_l$ onto the standard unit circle $\mathbb{S}=\partial\mathbb{D}$. By the Schwarz's reflection principle, we can extend $f$ such that the extension is analytic in a neighborhood of each $\gamma_l$ and hence $f'(w_{l,k})$ and $f'(w_{0,0})$ exist in Theorem \ref{thm:parameter} and we have $f(z)=1+f'(w_{l,k})(z-w_{l,k})+ \cdots$ near $w_{l,k}$ and $f(z)=1+f'(w_{0,0})(z-w_{0,0})+ \cdots$ near $w_{0,0}$. If we let $\mathbb{H}_{w_{l,k}}$ be the open half plane whose boundary is tangent to $\gamma_{l}$ at $w_{l,k}$ and it excludes $\gamma_l$ when $l \ge 1$ and includes $\gamma_0$ when $l=0$ and $k=0$, then since $f$ maps some neighborhoods of $w_{l,k}$ and $w_{0,0}$ in $\Omega$ into the same neighborhood of $1$ in $\mathbb{D}$, it follows that $\arg{f'(w_{l,k})}$ and $\arg{f'(w_{0,0})}$ are uniquely determined by $w_{l,k}$ and $w_{0,0}$ respectively. Theorem \ref{thm:parameter} shows, on the other hand, that the restrictions on $|f'(w_{l,k})|$ and $|f'(w_{0,0})|$ given by  (\ref{eqn:param}) will determine $f$ uniquely. \\

If we identify each $\gamma_l$ with the standard unit circle $\mathbb{S}$, then Theorem \ref{thm:parameter} implies that for $n=g+1$, we can define a surjective $n$ to $1$ map $G:\Omega \times \mathbb{S}^{n} \to \mathcal{P}_{n}(\Omega)$ which sends $(p,w_{0,0},\ldots,w_{g,0})$  to the unique proper map $f: \Omega \to \mathbb{D}$ satisfying (i) and (ii). While for $n>g+1$, we let $\mathbb{R}_{+}$ be the positive real axis. If $\nu=(n_{0},n_{1},\ldots, n_{g}) \in \mathbb{N}^{g+1}$ is a given boundary degree where $n=\sum_{l=0}^g n_l$, then we can define a surjective $(n\cdot n_0 !\cdots n_g !)$-to-$1$ map $G:\Omega  \times \mathbb{S}^{n} \times \mathbb{R}_{+}^{n-g-1}\to \mathcal{P}_{\nu}(\Omega)$ which sends 
\[(p,w_{0,0},\ldots,w_{g,n_{g}-1},\lambda_{0,1},\ldots,\lambda_{0,n_{0}-1},\ldots,\lambda_{g,1},\ldots, \lambda_{g,n_{g}-1})\]  to the unique proper map $f: \Omega \to \mathbb{D}$ satisfying conditions (i), (ii) and (iii) in Theorem \ref{thm:parameter}. We will further show that this map is continuous. This is the following corollary.
\begin{corollary}
\label{cor:parameter}
Let $\Omega$ be a circular domain of connectivity $g+1$. Let $n\ge g+1$ and $\nu=(n_{0},n_{1},\ldots, n_{g}) \in \mathbb{N}^{g+1}$ such that $n=\sum_{l=0}^g n_l$. We endow $\mathcal{P}_{\nu}(\Omega)$ with the topology of local uniform convergence. Then there exists a continuous $(n\cdot n_0 !\cdots n_g !)$-to-$1$ map from $\Omega  \times \mathbb{S}^{n} \times \mathbb{R}_{+}^{n-g-1}$ onto  $\mathcal{P}_{\nu}(\Omega)$. If $p\in \Omega$ and $\mathcal{P}_{\nu, p}(\Omega)$ is the subset of $\mathcal{P}_{\nu}(\Omega)$ which contains those proper holomorphic maps with a zero at $p$, then there is a continuous $(n_0 !\cdots n_g !)$-to-$1$ map from $\mathbb{S}^{n} \times \mathbb{R}_{+}^{n-g-1}$ onto  $\mathcal{P}_{\nu,p}(\Omega)$.
\end{corollary}

We remark that in \cite{bell2008}, Bell and Kaleem provide a different characterization of proper holomorphic maps of degree $n>g+1$ in terms of the preimages of $1$ using Grunsky maps. \\

This paper will be organized as follows. In Section \ref{sect:skpf}, we will give a brief introduction to the Schottky-Klein prime function and give some of its properties that we will need in the following sections. In Section \ref{sect:formula}, we will prove Theorems \ref{thm:expression} and \ref{thm:expression2}. In Section \ref{sect:cball}, we will use the previous results to establish Theorem \ref{thm:cball} and Corollary \ref{cor:parameter}. In Section \ref{sect:parameter}, we will prove Theorem \ref{thm:parameter}. Finally, in the last section, we will indicate how to associate a finite Blaschke product to a proper map in $\mathcal{P}(\Omega)$ when $\Omega$ is certain circular domain.
\section{The Schottky-Klein prime function}
\label{sect:skpf}

Let $\Omega$ be a $(g+1)-$connected non-degenerate circular domain and let $\gamma_{0}$ be the unit circle and $\gamma_1, \ldots , \gamma_g$ be the interior boundary circles. We suppose that $\gamma_{l}$ has center $q_l$ and of radius $r_l$ for each $l=0,\ldots,g$. For each $l=0,\ldots,g$, define 
\[ \varphi_l (z) = q_l + \dfrac{r_l^2}{\overline{z}-\overline{q_l}} \]
the reflection with respect to $\gamma_l$ and for $j=1,\ldots, g$, define 
\[ \theta_j (z) = \varphi_j( \varphi_0 (z)) =
q_j + \dfrac{r_j^2z}{1-\overline{q_j}z}.
\label{eqn:2refl}\]
As the composition of two reflections, $\theta_{j}$ are  M\"{o}bius transformations. The group $\Theta$ freely generated by $\theta_1 , \ldots , \theta_g$ is called the \textit{Schottky group} of $\Omega$. The domain $\Omega^{*}=\Omega \cup \gamma_0 \cup \varphi_0(\Omega)$ is a fundamental region associated with $\Theta$ called the \textit{Schottky double}. Then $\overline{\Omega^{*}}/\Theta$ is a compact Riemann surface of genus $g$. Let $\sigma_j$ be any line joining  a point $z \in \gamma_j$ with its identified point $\theta^{-1}_j(z) \in \varphi_0^{-1}(\gamma_j)$. The $2g$ curves $\gamma_1,\ldots, \gamma_g, \sigma_1, \ldots , \sigma_g$ form a basis of the fundamental group of $\overline{\Omega^{*}}/\Theta$. Then there exist $g$ holomorphic functions $v_1, \ldots , v_g$ satisfying the functional properties
\begin{equation}
\label{v_prop_1}
    \int_{\gamma_i} d v_j = \delta_{ij}
\end{equation}
and 
\begin{equation}
\label{v_prop_2}
    \int_{\sigma_i} d v_j = \tau_{ij}
\end{equation}
for some purely imaginary constants $\tau_{ij}$ and the Kronecker delta $\delta_{ij}$. These functions $v_{1},\ldots,v_{g}$ are called the \textit{integrals of the first kind} of $\Omega$.

We have the following properties:
\[\tau_{ij}=\tau_{ji}\mbox{ for } i,j=1,\ldots,g\]
and
\begin{equation}
\label{tauij}
    \tau_{ij}
    =v_j (\theta_i (z)) -v_j(z) 
\end{equation}
for any $z \in \overline{\Omega^{*}}$, for any $i= 1, \ldots ,g,$ for any $j=1, \ldots ,g$. See \cite{baker} and  \cite{crowdy19}.

Also, it is demonstrated in \cite{crowdy07} that 
\begin{equation}\label{vjsym} v_j( \overline{z}^{-1})  = \overline{ v_j (z)} \end{equation}
and 
\begin{equation}
    \label{har}
2 i
\begin{pmatrix}
\mathrm{Im} \left( v_1(z) \right) \\ 
\vdots \\ 
\mathrm{Im} \left( v_g(z) \right)
\end{pmatrix}
=
\begin{pmatrix}
\tau_{11} & \cdots &\tau_{1g} \\ 
\vdots & \ddots & \vdots \\ 
\tau_{g1} & \cdots & \tau_{gg} 
\end{pmatrix}
\begin{pmatrix}
u_1(z) \\ 
\vdots \\ 
u_g(z)
\end{pmatrix}
\end{equation}
where, for $j=1,\ldots,g$, $u_{j}$ denotes the harmonic measure of $\gamma_{j}$ (relative to $\Omega$).

For every $z,y \in \Omega$, the \textit{Schottky-Klein prime function} $\omega(z,y)$ is defined to be the unique function (of $z$) such that 
\begin{enumerate}
    \item $\omega(z,y)$ is analytic everywhere in $\Omega$;
    \item $\omega(z,y)$ has a single zero at each of the point $\{ \theta (y) \: : \: \theta \in \Theta \}$ (and is non-zero elsewhere);
    \item $\lim\limits_{z \to y} \dfrac{\omega(z,y)}{(z-y)} =1$, and,
    \item \begin{equation}\label{eq:omega_shift}\omega (\theta_j(z), y) = \exp \left( 2 \pi i \left( v_j(y) -v_j (z) \right) - \pi i \tau_{jj} \right) \sqrt{\dfrac{d \theta_j}{dz}}   \omega (z,y).\end{equation}
\end{enumerate}
The existence of such $\omega (z,y)$ is established in \cite{hejhal}. It would be useful to express $\omega (z,y)$ in the following infinite product form (see Chapter XII of \cite{baker}):
\begin{equation}
    \label{eq:def_skpf}
     \omega (z,y) = (z-y) \prod\limits_{\theta \in \Theta'} \dfrac{(z-\theta (y)) (y-\theta (z))}{(z-\theta (z))(y-\theta (y))}.
\end{equation}
Here $\Theta'$ is a maximal subset of $\Theta$ that excludes the identity element such that $\theta^{-1} \notin \Theta'$ whenever $\theta \in \Theta'$.
Note that the uniform convergence of this infinite product is still unknown for a general multiply connected planar
domain (see open problem 1 of \cite{crowdy2008geometric} and section 1 of \cite{mityushev2012}) and it is connected to the absolute convergence of the Poincar\'e theta series of the Schottky group (see \cite{baker}). It is known that for any circle decomposable Schottky group (see Section 4 of \cite{bobenko1989periodic} for its definition), its associated Poincar\'e theta series is absolutely convergent and hence the associated Schottky-Klein prime function $\omega(z,y)$ has the infinite product expression in (\ref{eq:def_skpf}). Examples of circle decomposable Schottky group are those Schottky groups with an invariant circle and in particular when all the centers $q_l$ are on the real axis (see page 54 of  \cite{burnside1891}). For the domain we are going to construct in the proof of Theorem \ref{thm:cball}, we will simply assume the centers of the boundaries of $\Omega$ are on the real axis even though the proof also works in the case when the radii $r_l$ are sufficiently less than the distances between the centers $|q_i-q_j|$ when $i\neq j$ (see page 55-58 of \cite{burnside1891}). 

Note that $\omega(z,y)$ is independent of the choice of $\Theta'$. In particular, the Schottky group $\Theta$ is the disjoint union of $\Theta'$, $(\Theta')^{-1}=\{\theta^{-1}:\theta\in\Theta'\}$ and the set containing the identity.

We will also need the following symmetry properties of $\omega$ (see \cite{crowdy19}): 
\begin{equation}
    \label{eq:omega_conju}
    \overline{\omega (z,y)} = - \overline{zy} \omega ( \overline{z}^{-1}, \overline{y}^{-1}).
\end{equation} 
From $(\ref{eq:def_skpf})$, we have \begin{equation}
    \label{eq:omega_exchan}
    \omega (z,y) = -\omega (y,z).
\end{equation} 
For more details about the Schottky-Klein prime function and its construction, see \cite{baker}, \cite{book_crowdy2020} and \cite{hejhal}.

Define 
\begin{equation}\eta(z,p)\equiv\frac{\omega(z,p)}{|p|\omega(z,\overline{p}^{-1})}\label{eq:circslitdisk} \end{equation}
and for any $z,p \in \Omega$ and each $l = 0, 1, \cdots, g$, define
\begin{equation}
    \label{eq:def_con}
    \eta_l (z,p) \equiv \sqrt{\dfrac{\varphi_l(p)}{p}} \dfrac{\omega (z,p)}{\omega (z, \varphi_l (p))}.
\end{equation}
In \cite{crowdy13}, it is proved that $\eta(\cdot,p)=\eta_0 (\cdot,p)$ is a conformal map from $\Omega$ to a circularly slit disk with zero $p$ mapping $\partial\mathbb{D}$ to $\partial \mathbb{D}$; and 
$\eta_j (\cdot,p)$ is a conformal map from $\Omega$ to a circularly slit disk with zero $p$ sending $\gamma_j$ to $\partial \mathbb{D}$. Note that such conformal maps are unique up to rotation. When we want to make the dependence on the domain $\Omega$ clear, we use the notation $\eta_{\Omega}$ for $\eta$ and $\eta_{\Omega,j}$ for $\eta_{j}$. For $j=1,\ldots,g$, we consider a M\"{o}bius transformation $\mu_{j}$ with $\mu_{j}(\gamma_{j})=\partial\mathbb{D}$ and $\mu_{j}(\Omega)\subset \mathbb{D}$. Then $\mu_{j}(\Omega)$ is also a non-degenerate circular domain and by the uniqueness of the conformal map onto a circular slit disk, \begin{equation}\eta_{\Omega,j}(\cdot,\cdot)=\lambda_{j}\eta_{\mu_{j}(\Omega)} (\mu_{j}(\cdot),\mu_{j}(\cdot))\label{eq:etajeta}\end{equation}
for some constant $\lambda_{j}$ where $|\lambda_{j}|=1$.

From the statements of Theorem \ref{thm:expression} and \ref{thm:expression2}, we can see that $\eta$ and $\eta_{j}$ act as the building blocks of proper holomorphic maps in $\mathcal{P}(\Omega)$. As mentioned in the introduction, this is analogous to the role of \[m(z,p)=\frac{z-p}{1-\overline{p}z}\] in finite Blaschke products (which are precisely the proper holomorphic maps in $\mathcal{P}(\mathbb{D})$). We will now prove (\ref{eqn:embedding}) which makes the link between the functions $m$ and $\eta$ explicit when the infinite product form in (\ref{eq:def_skpf}) holds.
\begin{proposition}\label{lem:embedding} Suppose that (\ref{eq:def_skpf}) holds, then
\[\eta(z,p)\equiv\prod_{\theta\in\Theta}\dfrac{m(\theta(z),p)}{m(\theta(1),p)}. \]
\end{proposition}
\begin{proof}
Fix $p\in\Omega\subset\mathbb{D}$.  Note that for any M\"{o}bius transformation $\mu$,
\[
\dfrac{\mu (z) - \mu (p)}{z-p} = \sqrt{\mu'(z)\mu'(p) }.
\]
In particular, for any $\theta\in\Theta$,
\begin{equation}
\label{eq:mob_shift}
    \dfrac{z-\theta(p)}{z-\theta(\overline{p}^{-1})}
    = \sqrt{\dfrac{\theta'(p)}{\theta'(\overline{p}^{-1})}}
    \left( \dfrac{\theta^{-1}(z)-p}{\theta^{-1}(z)-\overline{p}^{-1}}\right).
\end{equation}
Then from (\ref{eq:def_skpf}) and (\ref{eq:circslitdisk}), we have 
\begin{align*}
    \eta(z,p) &=
    \dfrac{\omega (z,p)}{|p| \omega (z,\overline{p}^{-1})}
    \\
    =&
    \left( \dfrac{1}{|p|} \right)
    \left( \dfrac{z-p}{z-\overline{p}^{-1}}\right)
    \prod_{\theta \in \Theta'}
     \dfrac{(z-\theta(p))(p-\theta(z))(z-\theta(z))(\overline{p}^{-1}-\theta(\overline{p}^{-1}))}{(z-\theta(\overline{p}^{-1}))(\overline{p}^{-1}-\theta(z))(z-\theta(z))(p-\theta(p))} \\
   =&
    \left( \dfrac{-1}{p} \right)
    \left( \dfrac{z-p}{ 1- \overline{p} z} \right)
    \prod_{\theta \in \Theta'}
 \dfrac{(z-\theta(p))(p-\theta(z))(\overline{p}^{-1}-\theta(\overline{p}^{-1}))}{(z-\theta(\overline{p}^{-1}))(\overline{p}^{-1}-\theta(z))(p-\theta(p))}.
\end{align*}
Using (\ref{eq:mob_shift}), the infinite product can be written as
\[ 
    \prod_{\theta \in \Theta'}
    \left(
   \overline{p} 
\sqrt{\dfrac{\theta'(p)}{\theta'(\overline{p}^{-1})}}
m (\theta^{-1} (z), p)
    \right)
    \left(
    \dfrac{\overline{p}(\overline{p}^{-1}-\theta(\overline{p}^{-1}))}{p-\theta(p)} m(\theta (z) ,p)
        \right)\]
and hence        
    \[
  \eta(z,p)=\left( \dfrac{-1}{p} m (z ,p) \right) 
    \prod_{\theta \in \Theta'}
    A_\theta m( \theta (z), p )m( \theta^{-1} (z), p )
\]
for some constants $A_\theta$ (depending on $p$ but not $z$). 

Now note that $\eta(1,p)=1$ and hence
\begin{align*}\eta(z,p)=&\dfrac{\eta(z,p)}{\eta(1,p)}\\=& \frac{m( z, p )}{m(1,p)}\dfrac{\prod\limits_{\theta \in \Theta'}
    A_\theta m( \theta (z), p )m( \theta^{-1} (z), p )}{\prod\limits_{\theta \in \Theta'}
    A_\theta m( \theta (1), p )m( \theta^{-1} (1), p )} \\ =& \frac{m( z, p )}{m(1,p)}\prod_{\theta \in \Theta'}
   \frac{m( \theta (z), p )}{m(\theta(1),p)}  \frac{m( \theta^{-1} (z), p )}{m(\theta^{-1}(1),p)}
   \\ =& \prod_{\theta\in\Theta}\dfrac{m(\theta(z),p)}{m(\theta(1),p)}
   \end{align*}
   where for the last equality, we have used the fact that
   \[\prod_{\theta \in \Theta'}
   \frac{m( \theta (z), p )}{m(\theta(1),p)} \mbox{ and } \prod_{\theta\in\Theta'}\frac{m( \theta^{-1} (z), p )}{m(\theta^{-1}(1),p)}\]
   both converge. This establishes the formula. 
\end{proof}
\section{Proving Theorems \ref{thm:expression} and \ref{thm:expression2}.}
\label{sect:formula}
As mentioned previously in the introduction, Wang and Yin proved in \cite{Wang} that there exists a proper holomorphic map $f: \Omega \to \mathbb{D}$ from $\Omega$ onto $\mathbb{D}$ of degree $n$ with zeros $p_{0},\ldots,p_{n}$ 
 if and only if the condition (\ref{eq:condition1}) is satisfied. The following lemma rewrites the condition (\ref{eq:condition1}) in terms of the integrals of the first kind $v_{1},\ldots,v_{g}$.
\begin{lemma}
\label{lem:con2}
Condition (\ref{eq:condition1}) is equivalent to 
\begin{equation}
    \label{eq:condition2}
         2 i\sum_{k=1}^n \mathrm{Im} \left(v_j (p_k) \right) = 
        \sum_{k=1}^n \tau_{jk} n_k 
     \qquad
     \mbox{for all $j=1,\ldots , g$.}
\end{equation}
\end{lemma}
\begin{proof}
Follows directly from (\ref{har}).
\end{proof}
For $z \in \Omega$ and $p_{0},\ldots,p_{n}\in\Omega$ satisfying condition (\ref{eq:condition1}), we consider the function given in Theorem \ref{thm:expression}:
\begin{equation}
    \label{eq:Fformula}
    f(z) =
\exp \left( -2 \pi i \sum\limits_{j=1}^{g} n_j v_j (z) \right) \prod_{k=1}^n \eta(z, p_k).
\end{equation}
The following lemmas give some properties of $f$ which we will require to prove Theorem \ref{thm:expression}.
\begin{lemma}
\label{lem2}
For any $z\in \partial\Omega$, $|f(z)|= 1$.
\end{lemma}

\begin{proof}
Suppose $z \in \gamma_0$. Then we have $z=\overline{z}^{-1}$ and by (\ref{har}), $\mathrm{Im}\left( v_j (z) \right)=0$ for all $j=1,\ldots, g$ and hence
\[\overline{\exp \left( -2 \pi i \sum\limits_{j=1}^{g} n_j v_j (z) \right)}
    \exp \left( -2 \pi i \sum\limits_{j=1}^{g} n_j v_j (z) \right)=1.\]
Also, it follows from (\ref{eq:omega_conju}) that
\[
 \overline{\eta(z, p)}\eta(z, p)=
 \dfrac{\overline{\omega (z, p)}\omega (z, p)}{|p|^2 \overline{\omega (z,  \overline{p}^{-1})}\omega (z,  \overline{p}^{-1}) }
= \dfrac{- \overline{z}p}{-\overline{z}\overline{p}^{-1}|p|^2}= 1.
\]
Using (\ref{eq:Fformula}) we deduce that, for $z\in\gamma_{0}$, $\overline{f(z)}f(z)=1$  and so $|f(z)|=1$.
Now, suppose that $z\in \gamma_i$ for some $i=1,\ldots, g$. Then we have $z=\theta_i (\overline{z}^{-1})$. It follows that 
\begin{align*}
    & \prod_{k=1}^n \dfrac{\overline{\omega (z, p_k)}\omega (z, p_k)}{|p_k|^2\overline{\omega (z,  \overline{p}_k^{-1})}\omega (z,  \overline{p}_k^{-1}) }
    \\
    = & \prod_{k=1}^n \dfrac{\overline{\omega (\theta_i (\overline{z}^{-1}), p_k)}\omega (z, p_k)}{|p_k|^2 \overline{\omega (\theta_i (\overline{z}^{-1}),  \overline{p}_k^{-1})}\omega (z,  \overline{p}_k^{-1}) }
    \\
    = & 
    \overline{\exp \left( 2\pi i \sum_{k=1}^n \left( v_i(p_k)-v_i(\overline{p_k}^{-1})\right) \right)} \prod_{k=1}^n   \dfrac{ \overline{\omega ( \overline{z}^{-1}, p_k)}\omega (z, p_k)}{|p_k|^2 \overline{\omega ( \overline{z}^{-1},  \overline{p}_k^{-1})}\omega (z,  \overline{p}_k^{-1}) }
    \tag{by (\ref{eq:omega_shift})}
        \\
    = & 
    \exp \left( 2\pi i \sum_{j=1}^g  \tau_{ij} n_j \right)
    \tag{by (\ref{vjsym}), (\ref{eq:condition2}) and (\ref{eq:omega_conju})}
\end{align*}
and 
\begin{align*}
    &
   \overline{ \exp \left( -2 \pi i \sum\limits_{j=1}^{g} n_j v_j (z) \right) }
   \exp \left( -2 \pi i \sum\limits_{j=1}^{g} n_j v_j (z) \right)
   \\
   =&
   \exp \left( 2 \pi i \sum\limits_{j=1}^{g} n_j \overline{v_j (\theta_i (\overline{z}^{-1}))} \right) 
   \exp \left( -2 \pi i \sum\limits_{j=1}^{g} n_j v_j (z) \right)
   \\
   =&
     \exp \left(2 \pi i \sum\limits_{j=1}^{g} n_j
     (v_j (z) - \tau_{ij})  
     - 2\pi i \sum\limits_{j=1}^{g} n_j
     v_j (z)  \right)  
     \tag{by (\ref{tauij}) and (\ref{vjsym})}
     \\
     =&
         \exp \left( -2\pi i \sum_{j=1}^g  \tau_{ij} n_j \right).
\end{align*}
Hence, $\overline{f(z)} f(z) = 1$ for $z \in \gamma_{i}$ for each $i=1, \ldots ,g$. Thus we have shown that $|f(z)|=1$ for all $z\in\partial \Omega$. 
\end{proof}

\begin{lemma}
\label{lem3}
For each $l=0, \cdots , g$, we have 
\[ \deg ( \left. f \right|_{\gamma_l}  )= n_l.\]
\end{lemma}

\begin{proof}
For each $k=1,\ldots ,n$, recall that $\eta (\cdot,p_k)$ is a conformal map from $\Omega$ to a circularly slit disk sending $\gamma_0$ to $\partial \mathbb{D}$ and $\gamma_j$ to a proper subarc $L_{jk}$ of a circle centred at $0$ for each $j=1,\ldots,g$. Notice that for if one travels along $\gamma_j$ once, the image of $\gamma_j$ under $\eta (\cdot,p_k)$ will move along $L_{jk}$ and then $-L_{jk}$.
For each $k=1,\ldots ,n$, considering $\eta=\eta (\cdot,p_k)$ as a change of variable, we have
\[ \dfrac{1}{2\pi i} \int_{\gamma_j}  d \log \eta (z,p_k)  = 
\dfrac{1}{2\pi i} \int_{L_{jk}-L_{jk}}  \dfrac{d \eta }{ \eta}
=
0\]
for all $j=1,\ldots,g$ and 
\[ \dfrac{1}{2\pi i} \int_{\gamma_0}  d \log \eta (z,p_k)  = 
\dfrac{1}{2\pi i} \int_{\partial \mathbb{D}}  \dfrac{d \eta }{ \eta}
=
1.\]
Using the argument principle, we have 
\begin{align*}
    \deg ( \left. f \right|_{\gamma_j}  )
    &=
    \dfrac{1}{2 \pi i} \int_{\gamma_j} d \log f(z)
    \\
    &= 
    \dfrac{1}{2 \pi i} \left(
    \sum\limits^g_{i=1} 2 \pi i n_i   \int_{\gamma_j} d v_i
    -
    \sum\limits^n_{k=1} \int_{\gamma_j}  d \log \eta (z,p_k)
    \right)
    \\
    &= n_j
\end{align*}
for each $j=1,\ldots,g$ and 
\begin{align*}
    \deg ( \left. f \right|_{\gamma_0}  )
    &=
    \dfrac{1}{2 \pi i} \int_{\gamma_0} d \log f(z)
    \\
    &=
    \dfrac{1}{2 \pi i} \int_{\partial \Omega} d \log f(z)
    - \sum\limits_{j=1}^g  \dfrac{1}{2 \pi i} \int_{\gamma_j} d \log f(z)
    \\
    &= n- \sum\limits_{j=1}^g n_j
    \\
    &= n_0.
\end{align*}
\end{proof}
We now prove Theorem \ref{thm:expression}.
\begin{proof}[Proof of Theorem \ref{thm:expression}]
By construction, it is clear that 
\[ f^{-1}(0) = \lbrace p_1 , \ldots , p_n \rbrace. \]
Lemma \ref{lem2} and the maximum modulus principle imply that 
$|f(z)|< 1$ for all $z\in\Omega$. This implies that $f$ is a finite map and hence is proper (see, for example, Section 9.3 of \cite{book_finite}). The condition of the degree of each boundary is implied by Lemma \ref{lem3}. 

As in \cite{Wang}, we obtain uniqueness as follows:  Suppose that $f_1$ is a proper holomorphic maps in $\mathcal{P}(\Omega)$ with zeros $p_1 , \ldots, p_n $. Then using the Schwarz reflection principle, we obtain holomorphic maps $\widetilde{f}_1$, $\widetilde{f}$ from the associated compact Riemann surface $\overline{\Omega^{*}}/ \Theta$ onto the Riemann sphere with zeros at $p_1 , \ldots, p_n $ and poles at $ \overline{p_{1}}^{-1} ,  \ldots, \overline{p_{n}}^{-1}$. Then $\dfrac{\widetilde{f}_1}{\widetilde{f}}$ has no zeros and poles and hence is constant. This constant   is equal to $e^{i\theta}$ for some real $\theta$ because $|f|=1$ on $\partial \Omega$ and $|f_1| \to 1$ as $z \to \partial \Omega$ (since $f_1$ is proper). This implies that $f_{1}=e^{i\theta}f$.
\end{proof}
To prove Theorem \ref{thm:expression2}, we will first need the following lemma which, in addition to (\ref{eq:etajeta}), gives a link between $\eta_{j}$ and $\eta$. 

\begin{lemma}
\label{lem:swapboundary}
For $\eta_{j}$ as defined in the statement of Theorem \ref{thm:expression2}, we have 
\[k(p)\eta_{j}(z,p)\equiv\exp(-2\pi i v_{j}(z)) \eta(z,p)\]
for some constant $k(p)$ not depending on $z$.
\end{lemma}
\begin{proof}
Since $\theta_{j}=\varphi_{j}\circ \varphi_{0}$, we have
\[ \overline{p}^{-1}=\theta_{j}^{-1}\circ \varphi_{j}(p).\]
Hence\begin{align*}
    \exp(-2\pi i v_{j}(z)) \eta(z,p)&= \exp(-2\pi i v_{j}(z))\frac{\omega(z,p)}{|p|\omega(z,\overline{p}^{-1})} \\
    &=\exp(-2\pi i v_{j}(z))\frac{\omega(z,p)}{|p|\omega(z,\theta_{j}^{-1}\circ \varphi_{j}(p) )} \\
    &= \widetilde{k}(p)\frac{\omega(z,p)}{\omega(z,\varphi_{j}(p))}. \tag{by (\ref{eq:omega_shift})}
\end{align*}
The result then follows from $(\ref{eq:def_con})$.
\end{proof}
\begin{proof}[Proof of Theorem \ref{thm:expression2}]
The expression for the proper holomorphic map $f$ follows directly from Lemma \ref{lem:swapboundary} and Theorem \ref{thm:expression}. It remains to show that if $\left\lbrace p_{l,k} : l=0,\ldots,g
\: \mbox{and} \:
 k=0,\ldots,n_{l}-1\right\rbrace$
satisfy condition (\ref{eq:condition3}) 
then there exists a proper holomorphic map with these points as zeros. 

Suppose that (\ref{eq:condition3}) is satisfied and let
\[
f(z)= R\prod_{l=0}^{g}\prod_{k=0}^{n_{l}-1} \eta_{l} (z,p_{l,k})
\]
where 
\[R= \left[\prod_{l=0}^{g}\prod_{k=0}^{n_{l}-1}\rho_{l,i}(p_{l,k})]^{-1}\right] \]
which does not depend on $i$. Then for $l=0,\ldots,g$ and any $w\in\gamma_{i}$, $|\eta_{l}(w,p_{l,k})|=\rho_{l,i}(p_{l,k})$ which implies that 
\[|f(w)|=\left|R\prod_{l=0}^{g}\prod_{k=0}^{n_{l}-1} \eta_{l} (w,p_{l,k})\right|=1.\]
So $|f(w)|=1$ for all $w\in\Omega$. As in the proof of Theorem \ref{thm:expression}, this implies that $f$ is proper. Then Theorem \ref{thm:expression} and Lemma \ref{lem:swapboundary} imply that $ \deg ( \left. f \right|_{\gamma_l}  )=n_{l} $ for $l=0,\ldots,g$. 
\end{proof}
\section{Disconnectedness of Carath\'{e}odory balls.}\label{sect:cball}
Suppose that $\Omega$ is a circular domain of connectivity $g+1$ ($g\ge 1$). In the introduction, we know from Grunsky's work (see \cite{grunsky1,grunsky2}) that the supremum in the definition of the M\"obius distance $c^{*}_{\Omega}(z_{0},z)$ is attained by some degree $g+1$ proper map in $\mathcal{P}_{g+1}(\Omega)$. 
Hence we can write 
\[c^{*}_{\Omega}(z_{0},z)=\sup\{|f(z)|: f\in\mathcal{P}_{g+1}(\Omega),f(z_{0})=0\}.\]
Moreover, by the uniqueness of proper maps up to rotation (cf. \cite{grunsky1950nachtrag}), we have
\[c^{*}_{\Omega}(z_{0},z)=\sup\{|f(z)|: f\in\mathcal{P}_{g+1}(\Omega),f(z_{0})=0,f(1)=1\}.\]

We will now construct the examples that will prove Theorem \ref{thm:cball}. 
Notice that since $c_{\Omega}=\tanh^{-1} c^{*}_{\Omega}$, we only need to establish (i) and (ii) in Theorem \ref{thm:cball} for the $c^{*}_{\Omega}$ balls $B_{c^*_{\Omega}}(z_{0};r_{1})$ and  $B_{c^*_{\Omega}}(z_{0};r_{2})$. Let $\widetilde{p}\in\Omega$ and \[\mathcal{Z}_{\Omega,\widetilde{p}}=\left\{\{p_{1},\ldots,p_{g}\} \subset \Omega :  u_{j}(\widetilde{p})+\sum_{k=1}^{g} u_j (p_k) = 1 
     \mbox{ for all $j=1,\ldots , g$}\right\}.\]
Thus, by Theorem \ref{thm:expression}, for each $P=\{p_{1},\ldots,p_{g}\}\in\mathcal{Z}_{\Omega,\widetilde{p}}$, there exists a unique degree $g+1$ proper holomorphic map of $\Omega$ onto $\mathbb{D}$ with zeros $\widetilde{p},p_{1},\ldots,p_{g} \in \Omega$ and with $1$ mapping to $1$. Denote this unique proper map by $\Phi_{\Omega}(\cdot;\widetilde{p}, P)$.
For any $\zeta \in \Omega$, there exists a unique $P_{\zeta}\in\mathcal{Z}_{\Omega,\widetilde{p}}$ such that $\Phi(\cdot,\widetilde{p},P_{\zeta})$ attains the supremum in the definition of $c^*_{\Omega}(\widetilde{p},\zeta)$ i.e.
\begin{equation}c^*_{\Omega}(\widetilde{p},\zeta)=|\Phi(\zeta;\widetilde{p},P_{\zeta})|.\label{eqn:extr}\end{equation}
Consider the $c^{*}_{\Omega}$ ball centered at $\widetilde{p}$ with radius $r$. This is 
\[B_{c^{*}_{\Omega}}(\widetilde{p};r)=\{z\in\Omega:c^*_{\Omega}(\widetilde{p},z)<r\}.\]
From the definition of the M\"obius distance, we  can write this as
\begin{align}
B_{c^{*}_{\Omega}}(\widetilde{p};r)&=\{z\in\Omega: \sup_{P\in\mathcal{Z}_{\Omega,\widetilde{p}}} |\Phi_{\Omega}(z;\widetilde{p},P)| <r\}\nonumber\\
&= \bigcap_{P\in\mathcal{Z}_{\Omega},\widetilde{p}} B_{P}(\widetilde{p};r) \label{eq:intersect}
\end{align}
where for each $P\in\mathcal{Z}_{\Omega,\widetilde{p}}$,
\[B_{P}(\widetilde{p};r):=\{z\in\Omega:  |\Phi_{\Omega}(z;\widetilde{p},P)| <r\}.\]

We now describe the general strategy used to prove Theorem \ref{thm:cball}. Theorem \ref{thm:expression2} implies that $\Phi(\cdot,\widetilde{p},P_{\zeta})$ can be written as an infinite product of the functions $\eta_{l}$. We will  find a domain $\Omega$ and a point $\zeta\in\Omega$ that is outside a neighbourhood of $\widetilde{p}$ so that one of the terms in the infinite product is arbitrarily small at the point $\zeta$ while the other terms are bounded. This will imply that $|\Phi(\zeta,\widetilde{p},P_{\zeta})|$ can also be made arbitrarily small. Using the fact that $\Phi(\widetilde{p},\widetilde{p},P_{\zeta})=0$, we will be able to find an $r>0$ so that $B_{P}(\widetilde{p};r)$ is disconnected. This, in turn, will allow us to establish Theorem \ref{thm:cball}.
\ \\

We use the notation from previous sections: 
$\Omega$ is a circular domain with boundary circles $\gamma_{0}=\partial\mathbb{D}$ and $\gamma_{1},\ldots, \gamma_{g} \subset \mathbb{D}$. In order to apply Proposition \ref{lem:embedding} in the following arguments, we will further assume that all the centers of $\gamma_{l}$ are on the real axis so that (\ref{eq:def_skpf}) holds.

For $q\in\Omega $, let $\eta_{\Omega}(\cdot,q)$ be a conformal map of $\Omega$ onto a circularly slit disk with $\partial\mathbb{D}$ mapping to the $\partial\mathbb{D}$ and $q$ mapping to $0$ and, for each $j=1,\ldots,g$ and $q\in\Omega$, let $\eta_{\Omega,j}(\cdot,q)$ be a conformal map of $\Omega$ onto a circularly slit disk with $\gamma_{j}$ mapping to the $\partial\mathbb{D}$ and $q$ mapping to $0$. We need the following lemma.

\begin{lemma} \label{lem:calc3}
For any given $j=1,\ldots,g$, and sequence $\{q_{n}\} \subset \Omega$ with $q_{n}\rightarrow \gamma_{j}$, we have 
\[|\eta_{\Omega,j}(\cdot,q_{n})|\rightarrow 1 \text{ as } n\rightarrow \infty\]
locally uniformly in $\Omega$. 
\end{lemma}
\begin{proof}
Suppose that $\{p_{n}\}$ is a sequence in $\Omega$ with $p_{n}\rightarrow \partial\mathbb{D}$. Using Proposition \ref{lem:embedding},
\[\eta_{\Omega}(w,p_{n})=\prod_{\theta \in \Theta}
   \frac{m( \theta (w), p_{n} )}{m(\theta(1),p_{n})}.\] 
Since $|m(z,p)|=1$ for any $p\in\partial\mathbb{D}$, we have for any compact set $K\subset\mathbb{D}$, $|m(\cdot,p_{n})|\rightarrow 1$ as $n\rightarrow \infty$ uniformly on $K$. Here, we note that for $n$ sufficiently large, $m(\theta(\cdot),p_{n})$ is bounded and does not vanish on $K$. Hence for all $\theta\in\Theta$,
\[\lim_{n\rightarrow\infty} \sup_{w\in K} \log |m(\theta(w),p_{n})| = 0.\]
Hence 
$\log |m(\theta(\cdot),p_{n})|$ is bounded on $K$ and the bounded convergence theorem then implies that 
\[ \sup_{w\in K} \lim_{n\rightarrow\infty} \sum_{\theta\in\Theta}\log |m(\theta(w),p_{n})| =  \sup_{w\in K} \sum_{\theta\in\Theta}\lim_{n\rightarrow\infty}\log |m(\theta(w),p_{n})|=0.\]
Hence
\[|\eta_{\Omega}(\cdot,p_{n})|\rightarrow 1 \]
converges uniformly on $K$ as $n \rightarrow\infty$. Then the result follows by setting $p_{n}=\mu_{j}(q_{n})$ and using (\ref{eq:etajeta}).
\end{proof}

We can then give a sufficient condition on $\Omega$ that would guarantee that $\Omega$ satisfies Theorem \ref{thm:cball}. This is contained in the following result.
\begin{theorem}\label{thm:suffcond}
Assume that $g \ge 2$. Suppose that for some distinct $j_{1},j_{2} \in \{1,\ldots,g\}$, there exist $\zeta\in \Omega$ satisfying
\begin{align}
\lim_{q\rightarrow \gamma_{j_{2}}} \left|\frac{ \eta_{\Omega}(\zeta,q)}{\eta_{\Omega}(w,q)}\right|&<\beta_{1}  \label{eqn:suffcond1}
\\
\lim_{q\rightarrow \partial\mathbb{D}} \left| \frac{\eta_{\Omega,j_{2}}(\zeta,q) }{\eta_{\Omega,j_{2}}(w,q)} \right|&< \beta_{2}
 \label{eqn:suffcond2} \end{align}
for some $0<\beta_{1}<1$, $0<\beta_{2}<\beta_{1}^{-1}$ and $w\in\gamma_{j_{1}}$.Then for some $\widetilde{p}\in\Omega$ sufficiently close to $\gamma_{j_{1}}$ and $r_{1}>r_{2}=c_{\Omega}(\widetilde{p},\zeta)$, 
\begin{enumerate}
\item[i)] $B_{c^*_{\Omega}}(\widetilde{p};r_{1})$ is disconnected and relatively compact in $\Omega$; 
\item[ii)]  $ \overline{B_{c^*_{\Omega}}(\widetilde{p};r_{2})} \subsetneq \overline{B}_{c^*_{\Omega}}(\widetilde{p};r_{2})$. 
\end{enumerate}
\end{theorem}
Note that the conditions in the statement of the theorem do not depend on the choice of $w\in\gamma_{j_{1}}$ since the functions $\eta_{\Omega}(\cdot,q)$ and $\eta_{\Omega,j}(\cdot,q)$ have constant modulus on each boundary component. Theorem \ref{thm:cball} will then follow from this theorem by finding a domain $\Omega$ in which (\ref{eqn:suffcond1}) and (\ref{eqn:suffcond2}) holds for some $\zeta\in\Omega$. 

\begin{proof}[Proof of Theorem \ref{thm:suffcond}]
Without loss of generality, we suppose $j_{1}=1$ and $j_{2}=2$.  As $\widetilde{p} \rightarrow \gamma_{1}$, we have $u_{1}(\widetilde{p})\rightarrow 1$ and $u_{j}(\widetilde{p})\rightarrow 0$ for each $j=2,\ldots,g$. Then as $\widetilde{p}\rightarrow \gamma_{1}$, the condition 
\[u_{j}(\widetilde{p})+\sum_{k=1}^{g} u_j (p_k) = 1  \mbox{ for all $j=1,\ldots , g$}\]
ensures that $u_{1}(p_{k})\rightarrow 0$ for each $k=1,\ldots ,g$. Hence, for each $P\in\mathcal{Z}_{\Omega,\widetilde{p}}$, every element of $P$ will be close to each of the other boundary components of $\Omega$. In particular, for each $\zeta\in\Omega$,  we can reorder the elements in $P_{\zeta}=\{p_{1,\zeta},\ldots, p_{g,\zeta}\} \in \mathcal{Z}_{\Omega,\widetilde{p}}$ such that as $\widetilde{p}\rightarrow \gamma_{1}$, we have
\begin{equation}p_{1,\zeta}\rightarrow \partial\mathbb{D} \mbox{ and } p_{j,\zeta}\rightarrow \gamma_{j} \mbox{ for } j=2,\ldots,g.\label{eqn:limbound}\end{equation}
By Theorem \ref{thm:expression2}, for $\zeta,w\in\Omega$ we can write
\[\left|\frac{\Phi_{\Omega}(\zeta;\widetilde{p},P_{\zeta})}{\Phi_{\Omega}(w;\widetilde{p},P_{\zeta})}\right|=\left|\frac{\eta_{\Omega}(\zeta,p_{2,\zeta})}{\eta_{\Omega}(w,p_{2,\zeta})}\right| \left|\frac{ \eta_{\Omega,1}(\zeta,\widetilde{p})}{ \eta_{\Omega,1}(w,\widetilde{p})}\right| \left|\frac{\eta_{\Omega,2}(\zeta,p_{1,\zeta})}{\eta_{\Omega,2}(w,p_{1,\zeta})}\right|  \prod_{j=3}^{g} \left| \frac{ \eta_{\Omega,j}(\zeta, p_{j,\zeta})}{ \eta_{\Omega,j}(w, p_{j,\zeta})}\right|\]
where the product $\prod_{j=3}^{g} \left| \frac{ \eta_{\Omega,j}(\zeta, p_{j,\zeta})}{ \eta_{\Omega,j}(w, p_{j,\zeta})}\right|$ is assumed to be $1$ when $g=2$.

By conditions (\ref{eqn:suffcond1}) and (\ref{eqn:suffcond2}), we can find $\widetilde{p}$ sufficiently close to  $\gamma_{1}$ such that for any $w\in \gamma_{1}$,
\begin{align*}|\eta_{\Omega}(\zeta,p_{2,\zeta})|&<\beta_{1} |\eta_{\Omega}(w,p_{2,\zeta})| 
\\ |\eta_{\Omega,2}(\zeta,p_{1,\zeta})|&<\beta_{2} |\eta_{\Omega,2}(w,p_{1,\zeta})|.\end{align*}
Hence 
\begin{align*}|\eta_{\Omega}(\zeta,p_{2,\zeta})||\eta_{\Omega,2}(\zeta,p_{1,\zeta})|&<  \beta_{1}\beta_{2}|\eta_{\Omega}(w,p_{2,\zeta})||\eta_{\Omega,2}(w,p_{1,\zeta})|\\&<|\eta_{\Omega}(w,p_{2,\zeta})||\eta_{\Omega,2}(w,p_{1,\zeta})|\end{align*}
since $\beta_{1}\beta_{2}<1$.

Taking such a value of $\widetilde{p}$ and $0<\epsilon<1$ sufficiently close to $1$ such that $\beta_{1}\beta_{2}<\epsilon$, define $\rho>0$ by
\[\rho= \epsilon |\eta_{\Omega}(w,p_{2,\zeta})||\eta_{\Omega,2}(w,p_{1,\zeta})|\]for some $w\in\gamma_{1}$
such that 
\[|\eta_{\Omega}(\zeta,p_{2,\zeta})||\eta_{\Omega,2}(\zeta,p_{1,\zeta})|<\rho< |\eta_{\Omega}(\widetilde{p},p_{2,\zeta})||\eta_{\Omega,2}(\widetilde{p},p_{1,\zeta})|.\]
For the reasons mentioned previously, the definition of $\rho$ does not depend on the choice of $w\in\gamma_{1}$. Let
\[L=\{u\in\Omega: |\eta_{\Omega}(u,p_{2,\zeta})||\eta_{\Omega,2}(u,p_{1,\zeta})|= \rho\}.\] 
Then $L$ is a closed (and hence compact) subset of $\Omega$ and any path from $\widetilde{p}$ to $\zeta$ in $\Omega$ must intersect $L$. Notice that we have for any $u\in L$,
\begin{equation} |\eta_{\Omega}(\zeta,p_{2,\zeta})||\eta_{\Omega,2}(\zeta,p_{1,\zeta})|<  \frac{\beta_{1}\beta_{2}}{\epsilon} |\eta_{\Omega}(u,p_{2,\zeta})||\eta_{\Omega,2}(u,p_{1,\zeta})|.\label{eqn:suffcond3} \end{equation}
Then (\ref{eqn:suffcond3}) implies that we can find $\widetilde{p}$ sufficiently close to $\gamma_{1}$ such that for any $u\in L$, 
 \begin{equation}
\left|\frac{\eta_{\Omega}(\zeta,p_{2,\zeta})}{\eta_{\Omega}(u,p_{2,\zeta})}\right|  \left|\frac{\eta_{\Omega,2}(\zeta,p_{1,\zeta})}{\eta_{\Omega,2}(u,p_{1,\zeta})}\right| <\frac{\beta_{1}\beta_{2}}{\epsilon}<1.
 \label{eqn:suffcond3a} \end{equation}
 From (\ref{eqn:limbound}) and Lemma \ref{lem:calc3}, for each $j=3,\ldots,g$,
\[\left|\frac{ \eta_{\Omega,1}(\zeta,\widetilde{p})}{ \eta_{\Omega,1}(u,\widetilde{p})}\right|\rightarrow 1 \mbox{ and }  \left| \frac{ \eta_{\Omega,j}(\zeta, p_{j,\zeta})}{ \eta_{\Omega,j}(u, p_{j,\zeta})}\right|\rightarrow 1 \mbox{ as } \widetilde{p}\rightarrow \gamma_{1}\]
and this convergence is uniform for $u$ in the compact set $L$. Hence, for $u\in L$, 
\[\lim_{\widetilde{p}\rightarrow \gamma_{1}}  \left|\frac{\Phi_{\Omega}(\zeta;\widetilde{p},P_{\zeta})}{\Phi_{\Omega}(u;\widetilde{p},P_{\zeta})}\right|= \lim_{\widetilde{p}\rightarrow \gamma_{1}} \left|\frac{\eta_{\Omega}(\zeta,p_{2,\zeta})}{\eta_{\Omega}(u,p_{2,\zeta})}\right|  \left|\frac{\eta_{\Omega,2}(\zeta,p_{1,\zeta})}{\eta_{\Omega,2}(u,p_{1,\zeta})}\right|. \]
 It follows from (\ref{eqn:suffcond3a}) that for $\widetilde{p}$ sufficiently close to $\gamma_{1}$,
\[\left|\frac{\Phi_{\Omega}(\zeta;\widetilde{p},P_{\zeta})}{\Phi_{\Omega}(u;\widetilde{p},P_{\zeta})}\right| <1\]
for all $u\in L$.
Since $L$ is compact, we can find $0<r_{1}<1$ such that
 \[|\Phi_{\Omega}(\zeta;\widetilde{p},P_{\zeta})|<r_1<|\Phi_{\Omega}(u;\widetilde{p},P_{\zeta})|\]
 for all $u\in L$ .
Notice that (\ref{eqn:extr}) implies that
$c^*_{\Omega}(\widetilde{p},\zeta)<r_{1}$
 and therefore $\widetilde{p},\zeta \in B_{c^*_{\Omega}}(\widetilde{p};r_{1})$. As $r_1<|\Phi_{\Omega}(u;\widetilde{p},P_{\zeta})|$ for all $u\in L$, $L\cap  B_{P_{\zeta}}(\widetilde{p};r_{1}) = \emptyset$. Then (\ref{eq:intersect}) implies that $L\cap B_{c^*_{\Omega}}(\widetilde{p};r_{1})=\emptyset$. Since any path in $\Omega$ from $\widetilde{p}$ to $\zeta$ must intersect $L$, we conclude that $B_{c^*_{\Omega}}(\widetilde{p};r_{1})$ is disconnected.
 Also $B_{c^*_{\Omega}}(\widetilde{p};r_{1})$ is relatively compact as its closure is bounded.
 
 Finally, since $B_{c^*_{\Omega}}(\widetilde{p};r_{1})$ is disconnected, there exists a connected component of it which does not contain $\widetilde{p}$. Let $S$ be such a component and $C$ be the union of all the other components so that $B_{c^*_{\Omega}}(\widetilde{p};r_{1})=S\cup C$ and
$\overline{S} \cap C =\emptyset$. Let $r_{2}= \inf\{c^*_{\Omega}(\widetilde{p},z): z \in \overline{S}\}$. Notice that $\overline{B_{c^*_{\Omega}}(\widetilde{p};r_{1})}$ is compact and hence so is $\overline{S}$. As $c^*_{\Omega}(\widetilde{p},z)$ is a continuous function in $z$, there exists some $\xi \in \overline{S}$ such that $0<r_2=c^*_{\Omega}(\widetilde{p},\xi)<r_1$. Clearly, $\xi \in \overline{B}_{c^*_{\Omega}}(\widetilde{p};r_{2})$, $\overline{S}\cap B_{c^*_{\Omega}}(\widetilde{p};r_{2})=\emptyset$ and $\overline{B_{c^*_{\Omega}}(\widetilde{p};r_{2})} \subset B_{c^*_{\Omega}}(\widetilde{p};r_{1})=S\cup C$. It follows from  $\overline{S}\cap B_{c^*_{\Omega}}(\widetilde{p};r_{2})=\emptyset$ that $\overline{B_{c^*_{\Omega}}(\widetilde{p};r_{2})}\cap S=\emptyset$ and hence $\overline{B_{c^*_{\Omega}}(\widetilde{p};r_{2})}\subset C$. Since 
$\overline{S} \cap C =\emptyset$, we have $\overline{B_{c^*_{\Omega}}(\widetilde{p};r_{2})}\cap \overline{S}=\emptyset$.
As $\xi \in \overline{S}$, we must have  $\overline{B_{c^*_{\Omega}}(\widetilde{p};r_{2})}\subsetneq \overline{B}_{c^*_{\Omega}}(\widetilde{p};r_{2})$.
 
 \end{proof}
  It remains to find a domain $\Omega$ and a point $\zeta\in\Omega$ that satisfies (\ref{eqn:suffcond1}) and (\ref{eqn:suffcond2}). We first need to introduce further notation. For $\widehat{r}>0$, let $\{\widetilde{\gamma}_{r}\subset \mathbb{D}:r\in(0,\widehat{r}]\}$ be a family of concentric circles in $\mathbb{D}$ such that $\widetilde{\gamma}_{r}$ has radius $r$. 
 Let $\gamma_{1},\ldots ,\gamma_{g-1} \subset \mathbb{D}$ be disjoint circles so that they together with the unit circle bound a circular domain $\Omega^*$. If we further assume that $\gamma_{j} \cap \widetilde{\gamma}_{r}=\emptyset$ for any $r\in(0,\widehat{r}]$ and $j=1,\ldots, g-1$, then we define the family of $g+1$ connected domains $\{\Omega_{r}:r\in(0,\widehat{r}]\}$ such that $\Omega_{r}$ is the circular domain bounded by the unit circle and the circles $\gamma_{1},\ldots ,\gamma_{g-1}$ and $\widetilde{\gamma}_{r}$.   We first need the following lemma.
 \begin{lemma} \label{lem:calc} 
 Suppose that $r\in(0,\widehat{r}]$ and let $\eta_{r}=\eta_{\Omega_{r}}$. Then as $r\rightarrow 0$,
 \[\eta_{r}(\cdot,\cdot)\rightarrow \eta_{\Omega^{*}}(\cdot,\cdot)\]
 uniformly in $\overline{\Omega_{\widehat{r}}}\times\overline{\Omega_{\widehat{r}}}$. 
 \end{lemma}
 \begin{proof}
For $r\in(0,\widehat{r}]$, let $\Theta_{r}$ denote the Schottky group associated with $\Omega_{r}$. Similarly, let $\Theta^{*}$ denote the Schottky group associated with $\Omega^{*}$. Then (using the notation from Section \ref{sect:skpf}),  $\Theta^{*}=\langle \theta_{1},\ldots ,\theta_{g-1},\theta_{1}^{-1},\ldots ,\theta_{g-1}^{-1} \rangle$ and $\Theta_{r}= \langle \theta_{1},\ldots ,\theta_{g-1}, \widetilde{\theta}_{r},\theta_{1}^{-1},\ldots ,\theta_{g-1}^{-1}, \widetilde{\theta}_{r}^{-1} \rangle$ where $\widetilde{\theta}_{r}$ is the element in $\Theta_{r}$ associated with $\widetilde{\gamma}_{r}$.   Using Proposition \ref{lem:embedding},
\begin{equation}\eta_{r}(w,p)=\eta_{\Omega^{*}}(w,p)\prod_{\theta\in T_{r}}  \frac{m(\theta(w),p)}{m(\theta(1),p)}\label{eqn:calc}\end{equation}
where $T_{r}=\Theta_{r}\setminus \Theta^{*}$. Notice that $\widetilde{\theta_{r}}(z)=\alpha + \dfrac{r^2z}{1-\overline{\alpha}z}$ where $\alpha\in \mathbb{D}\backslash \Omega_r$ is the center of $\widetilde{\gamma}_{r}$. It follows that \[\lim_{r\rightarrow 0}\sup_{w\in \overline{\Omega_{\widehat{r}}} } \left| \widetilde{\theta}_{r}(w)-\alpha \right| = 0 .\]
Similarly, 
\[\lim_{r\rightarrow 0}\sup_{w\in \overline{\Omega_{\widehat{r}}} } \left| \widetilde{\theta}_{r}^{-1}(w)-\overline{\alpha}^{-1}\right|= 0 .\]
 Together with the observation that any $\theta\in T_{r}$ has at least one $\widetilde{\theta}_{r}$ or  $\widetilde{\theta}_{r}^{-1}$ in its compositions of the generators of $\Theta_{r}$, we have for any $\theta\in T_{r}$, there exists some complex number $\nu$ such that 
 \[\lim_{r\rightarrow 0}\sup_{w\in \overline{\Omega_{\widehat{r}}} } \left| \theta(w)-\nu \right| = 0 \]
 and in particular $\lim_{r\rightarrow 0}\theta(1)=\nu$. Hence
 \[ \lim_{r\rightarrow 0}\sup_{(w,p)\in \overline{\Omega_{\widehat{r}}}^{2}}\log \left|\frac{m(\theta(w),p)}{m(\theta(1),p)}\right|=0.\]
Here we use the fact that $\dfrac{m(\theta(w),p)}{m(\theta(1),p)}$ does not vanish for $(w,p)\in \overline{\Omega_{\widehat{r}}}^{2}$ for $\theta\in T_{r}$. As 
 $\log \left|\dfrac{m(\theta(w),p)}{m(\theta(1),p)}\right|$ is bounded for
 $(w,p)\in \overline{\Omega_{\widehat{r}}}^{2}$, by the bounded convergence theorem, we deduce that
\[\sup_{(w,p)\in \overline{\Omega_{\widehat{r}}}^{2}} \lim_{r\to 0}\sum_{\theta \in T_{r}} \log\left| \frac{m(\theta(w),p)}{m(\theta(1)),p)} \right|=\sup_{(w,p)\in \overline{\Omega_{\widehat{r}}}^{2}} \sum_{\theta \in T_{r}} \lim_{r\to 0}\log\left| \frac{m(\theta(w),p)}{m(\theta(1)),p)} \right| = 0.
\]
Notice that $T_r$ is a countable set and we can order elements of $T_r$ as $\theta_1,\theta_2,...$ where each $\theta_i$ depends on $r$. Then $\sum_{\theta \in T_{r}}$ is understood to be  $\sum_{i=1}^{\infty}$.
Thus,
\[\sup_{(w,p)\in \overline{\Omega_{\widehat{r}}}^{2}} \left|\left[\prod_{\theta\in T_{r}}  \frac{m(\theta(w),p)}{m(\theta(1),p)}\right]-1 \right| \rightarrow 0 \text{ as } r\rightarrow 0.\]
Therefore, using (\ref{eqn:calc}),
\begin{align*}&\sup_{(w,p)\in \overline{\Omega_{\widehat{r}}}^{2}} \left|\eta_{r}(w,p)-\eta_{\Omega^{*}}(w,p)\right|\\ =& \sup_{(w,p)\in \overline{\Omega_{\widehat{r}}}^{2}} \left|\eta_{\Omega^{*}}(w,p)\right| \left|\left[\prod_{\theta\in T_{r}}  \frac{m(\theta(w),p)}{m(\theta(1),p)}\right]-1\right|
\\ \leq& \sup_{(w,p)\in \overline{\Omega_{\widehat{r}}}^{2}}  \left|\left[\prod_{\theta\in T_{r}}  \frac{m(\theta(w),p)}{m(\theta(1),p)}\right]-1\right|
\rightarrow 0\end{align*}
as $r\rightarrow 0$. This completes the proof.
\end{proof}
 We also need the following lemma.
\begin{lemma}\label{lem:calc2}
Let $\widetilde{\eta}_{r}$ denote the conformal map of $\Omega_{r}$ onto a circularly slit disk with $\widetilde{\gamma}_{r}$ mapping onto $\partial\mathbb{D}$ and 1 mapping to 1. Let $A\subset \Omega_{\widehat{r}}$ and let $K$ be a compact set in $\overline{\Omega_{\widehat{r}}}$ such that $K\cap \overline{A} = \emptyset$. Then there exists $M_{1},M_{2}>0$ such that
\[ \inf\{|\widetilde{\eta}_{r}(w, q)||:w\in K, q\in A, r\in(0,\widehat{r}]\} >M_{1}\]
and
\[ \sup\{|\widetilde{\eta}_{r}(w, q)||:w\in K, q\in A, r\in(0,\widehat{r}]\} <M_{2}\]
 \end{lemma}
 \begin{proof}
 Consider the family of functions 
 \[\mathcal{F}=\{\widetilde{\eta}_{r}(\cdot, q): q\in A, r\in(0,\widehat{r}] \}.\]
 This family is uniformly bounded on a neighbourhood of $\overline{\Omega_{\widehat{r}}}$ and is non-zero on $K$ as the zeros of functions in $\mathcal{F}$ are in $A$ which is disjoint from $K$. By Montel's theorem, $\mathcal{F}$ forms a normal family which means $\mathcal{F}$ is pre-compact in the topology of local uniform convergence and clearly $M_2$ exists.  Now if $h_0 \in \overline{\mathcal{F}}$ attains $\inf_{h\in\overline{\mathcal{F}}}\min\{|h(w)|:w \in K\}$.
 Then $h_0$ is the limit of some sequence in $\mathcal{F}$ and hence $h_0$ cannot have any zero in $K$. Thus $\min\{|h_0(w)|:w \in K\}>0$ and $M_1$ exists. 
\end{proof}
 We can now prove Theorem \ref{thm:cball}.
 \begin{proof}[Proof of Theorem \ref{thm:cball}] 
As mentioned previously, we need to show the existence of $\Omega$ and $\zeta$ that satisfies (\ref{eqn:suffcond1}) and (\ref{eqn:suffcond2}) in Theorem \ref{thm:suffcond}.  
Using the notation in Lemma \ref{lem:calc}, note that  
\[\lim_{\zeta, q\rightarrow \alpha} \left|\frac{\eta_{\Omega^{*}}(\zeta,q)}{\eta_{\Omega^{*}}(w,q)}\right|\rightarrow 0 \]
for $w\in\gamma_{1}$. Here we have used the the continuity of $\eta_{\Omega^{*}}(\cdot, \cdot)$ and the fact that $\eta_{\Omega^{*}}(q,q) = 0$.
Then Lemma \ref{lem:calc} implies that for $\epsilon>0$, we can find $r>0$ sufficiently small and $\zeta, \ q$ sufficiently close to $\widetilde{\gamma}_{r}$ such that
\[\left|\frac{ \eta_{r}(\zeta,q)}{\eta_{r}(w,q)}\right|<\epsilon\]
for $w\in\gamma_{1}$.
Now, Lemma \ref{lem:calc2} implies that 
\[\left| \frac{\widetilde{\eta_{r}}(\zeta,q) }{\widetilde{\eta}_{r}(w,q)} \right|\]
is bounded. Hence for $\epsilon$ sufficiently small and  for $\gamma_{j_{1}}=\gamma_{1}$ and $\gamma_{j_{2}}=\widetilde{\gamma}_{r}$, (\ref{eqn:suffcond1}) and (\ref{eqn:suffcond2}) hold for this $\Omega_{r}$ and $\zeta$. This completes the proof.
\end{proof} 
\section{The parameter space of proper holomorphic maps.}
\label{sect:parameter}
We will now prove Theorem \ref{thm:parameter}. Recall that, from the definition of harmonic measure, $u_{j}(w)=\delta_{jl}$ for $w\in \gamma_{l}$. Hence, if
$w_{l,k} \in\gamma_{l}$ for $l=0,\ldots,g$; $k=0,\ldots,n_{l}-1$, then for each $1\le j \le g$, we have
\begin{equation}\label{eq:conditionw}
\sum_{l=0}^{g} \sum_{k=0}^{n_{l}-1} u_{j}(w_{l,k})=\sum_{k=0}^{n_{j}-1} u_{j}(w_{j,k})=n_{j}.
\end{equation}
In other words, the condition (\ref{eq:condition1}) is satisfied when the points $p_i$ are replaced by the boundary points $w_{l,k}$. The idea of the proof is to use the implicit function theorem to find points in $\Omega$ sufficiently close to the boundary such that the condition (\ref{eq:condition1}) is satisfied by these points. We are then guaranteed a proper map with zeros at these points and with the required degree on each boundary component of $\Omega$. Once we obtain such a proper map, we can then use this to find a proper map which satisfies the desired properties. Uniqueness of the proper map will also follow as a consequence of the implicit function theorem. We will first need the following lemmas.
 \begin{lemma}\label{lem:allpoints}
 Suppose $f\in \mathcal{P}_{\nu}(\Omega)$ and $\zeta\in\mathbb{D}$. Let $q_{1},\ldots ,q_{n}\in\Omega$ be the preimages of $\zeta$ under $f$. Then $q_{1},\ldots, q_{n}$ satisfy the condition (\ref{eq:condition1}). Conversely for $\zeta\in\mathbb{D}$, suppose that $q_{1},\ldots, q_{n}\in\Omega$ satisfy the condition (\ref{eq:condition1}), then there exists $f\in \mathcal{P}_{\nu}(\Omega)$ with $f(1)=1$ such that the pre-image of $\zeta$ under $f$ is $q_{1},\ldots,q_{n}$.
 Consequently, for any $p\in\Omega$, there exists $f\in \mathcal{P}_{\nu}(\Omega)$ with $f(1)=1$, $f(p)=0$ and \[f(q_{1})=f(q_{2})=\cdots=f(q_{n}).\]
 \end{lemma}
 \begin{proof}
Let $\mu$ be a M\"{o}bius transformation that fixes $\mathbb{D}$ and maps $\zeta$ to $0$. Then $\mu\circ f\in\mathcal{P}_{\nu}(\Omega)$ and has zeros $q_{1},\ldots, q_{n}$. Hence, $q_{1},\ldots, q_{n}$ satisfy the  condition  (\ref{eq:condition1}).

Conversely, if $q_{1},\ldots, q_{n}\in\Omega$ satisfy the condition (\ref{eq:condition1}), then there exists a proper holomorphic map $g\in\mathcal{P}_{\nu}(\Omega)$ such that $g(q_{k})=0$ for $k=1,\ldots,n$. Let $\mu$ be a M\"{o}bius transformation that fixes $\mathbb{D}$ and maps 0 to $\zeta$ and $g(1)$ to $1$. Then $f=\mu\circ g$ satisfies the required conditions. For the final statement, it suffices to consider $\widetilde{\mu}\circ f$ where $\widetilde{\mu}$ is a M\"{o}bius transformation  that fixes $\mathbb{D}$ and maps $f(p)$ to 0 and $f(1)$ to $1$.
 \end{proof} 
 
 We shall need the following version of Hopf's lemma (see p.44 of \cite{garnett2005harmonic}).
 
 \begin{lemma}\label{lem:Hopf}
 Suppose $D$ is a non-degenerate  finitely connected domain in $\mathbb{C}$ whose boundary $\partial D$ consists of Jordan curves. Let $\gamma \subset \partial D$ be an analytic arc and $u$ be a harmonic function in $D$. If $\lim_{z\to \zeta}u(z)=0$ for all $\zeta \in \gamma$, then there is an open set $W$ containing $\gamma \cup D$ such that $u$ extends to be harmonic on $W$. If also $u(z)<0$ in $D$, then the 
 \[\frac{\partial u}{\partial n}(\zeta) >0 \]
 for all $\zeta \in \gamma$. Here $n$ is the unit normal vector pointing out from the domain $D$.
 \end{lemma}
 
 \begin{lemma}\label{lem:invertible}
Suppose that $w_{1},\ldots,w_{g}$ satisfy $w_{j}\in \gamma_{j}$. Let $U$ be the $g$-by-$g$ matrix whose entry in the $j$-th row and $k$-th column is $u_{j}'(w_{k})=\frac{\partial u_j}{\partial n}(w_{k})$. Then $U$ is non-singular.
\end{lemma}
\begin{proof}
We will show that the rows of $U$ are linearly independent over $\mathbb{R}$. The rows of $U$ are the vectors 
\[x_{j}=\left(u_{j}'(w_{1}),\ldots, u_{j}'(w_{g})\right) \quad \text{ for } j=1,\ldots,g.\]
Take $\alpha_{1},\ldots,\alpha_{g}\in\mathbb{R}$ and define
\[\Psi(z)=\sum_{j=1}^{g}\alpha_{j}u_{j}(z).\]
Then $\Psi$ is harmonic and takes the value $\alpha_{j}$ on $\gamma_{j}$ for $j=1,\ldots,g$ and $0$ on $\gamma_{0}$.  
Then 
\begin{align*}
& & & &  
\alpha_{1}x_{1}+\cdots + \alpha_{g}x_{g}&=0
&
\\
&\Longleftrightarrow  & & & \sum_{j=1}^{g} \alpha_{k} u_{j}'(w_{k}) &= 0  
&\quad \text{for all }k=1,\ldots, g \\
&\Longleftrightarrow & 
& & 
\frac{\partial \Psi}{\partial n}(w_{k}) &=0
&\quad \text{for all } k =1,\ldots, g. 
\end{align*}
By relabelling the boundary components of $\Omega$ if necessary, we can assume that $\alpha_{1}\geq \alpha_{j}$ for $j=2,\ldots,g$. Now, suppose (for a contradiction) that $\Psi$ is non-constant. The maximum principle of harmonic functions implies that $\Psi(z)<\alpha_{1}$ for all $z\in\Omega$ and hence, by Lemma \ref{lem:Hopf}, $\frac{\partial \Psi}{\partial n}(w)\neq 0$ for $w\in\gamma_{1}$. This is a contradiction since $\frac{\partial \Psi}{\partial n}(w_{1})=0$. Hence we must have $\Psi(z)$ is constant and thus $\Psi(z)\equiv 0$ (because $\Psi(z)\equiv 0$ on $\gamma_0$). Thus we conclude that $\alpha_1=\cdots=\alpha_{g}=0$ and so $x_{1},\ldots,x_{g}$ are linearly independent over $\mathbb{R}$.
\end{proof}

We now prove Theorem \ref{thm:parameter}. 
\begin{proof}[Proof of Theorem \ref{thm:parameter}]
Take $p\in\Omega$ and $w_{l,k}\in \gamma_{l}$ for $l=0,\ldots,g$ and $k=0,\ldots ,n_{l}-1$. 
We will construct a proper holomorphic map 
$F\in\mathcal{P}_{\nu}(\Omega)$ satisfying $F(p)=0$ and $F(w_{l,k})=1$ for all $l=0,\ldots,g$ and $k=0,\ldots ,n_{l}-1$.

For each $l=0,\ldots,g$ and $k=0,\ldots,n_{l}-1$, let ${\bf n}_{l,k}$ be the unit normal vector of $\gamma_l$ at $w_{l,k}$, pointing towards the interior of $\Omega$. Then for $r_{l,k} \in (-\varepsilon,\varepsilon)$, we consider the normal line of $\gamma_l$ at $w_{l,k}$, 
\[w_{l,k}+r_{l,k}{\bf n}_{l,k}.\]
We will apply the implicit function theorem to the following system of equations
\begin{equation}\left[\sum_{l=0}^{g} \left(u_{j}(w_{l,0}+r_{l,0}{\bf n}_{l,0})+\sum_{k=1}^{n_{l}-1} u_{j}(w_{l,k}+r_{l,k}{\bf n}_{l,k})\right)\right]-n_{j}=0,\label{eqn:implicitcondition}\end{equation}
where $j=1,\ldots,g$. Here and in the following argument, we assume that for any $l$ with $n_{l}=1$, $\sum_{k=1}^{n_{j}-1}$ is taken as an empty sum.
The left-hand side of this system of equations is well-defined for sufficiently small positive $\epsilon$ because by  Lemma \ref{lem:Hopf} or the Schwarz reflection principle, we can extend each $u_j$ to a harmonic function in a neighborhood $W$ of $\overline{\Omega}$. Notice that $u_j(z)>1$ for any $z \in W$ and inside the open disk bounded by $\gamma_j$ while $u_j(z)<0$ for any $z \in W$ and in the open disk bounded by $\gamma_i$  for any $i \neq 0,j$. Let \begin{align*}&F_j(r_{1,0},\ldots r_{g,0},r_{0,0},\ldots ,r_{0,n_{0}-1},r_{1,1}\ldots ,r_{1,n_{1}-1},\ldots ,r_{g,1}\ldots ,r_{g,n_{g}-1})\\ &=\left[\sum_{l=0}^{g} \left(u_{j}(w_{l,0}+r_{l,0}{\bf n}_{l,0})+\sum_{k=1}^{n_{l}-1} u_{j}(w_{l,k}+r_{l,k}{\bf n}_{l,k})\right)\right]-n_{j}\end{align*}
and $F: (-\varepsilon, \varepsilon)^g \times (-\varepsilon,\varepsilon)^{n-g} \subset \mathbb{R}^g \times \mathbb{R}^{n-g} \to \mathbb{R}^g$ be defined by $F=(F_1,\ldots ,F_g)$. From equation (\ref{eq:conditionw}), we have at $({\bf 0},{\bf 0}) \in  \mathbb{R}^g \times \mathbb{R}^{n-g}$, $F({\bf 0},{\bf 0})={\bf 0}$. Notice that for $1 \le i,j \le g$,  $\frac{\partial F_j}{\partial r_{i,0}}({\bf 0},{\bf 0})=\frac{\partial u_j}{\partial n}(w_{i,0})$ by the definition of normal derivative and hence by Lemma \ref{lem:invertible}, $\det((\frac{\partial F_j}{\partial r_{i,0}}({\bf 0},{\bf 0}))) \neq 0$. 
Thus the implicit function theorem guarantees that there exist $0<\varepsilon_0, \varepsilon_1 <\varepsilon$ and $g:(-\varepsilon_0, \varepsilon_0)^g \to (-\varepsilon_1, \varepsilon_1)^g$ of class
$\mathcal{C}^1$ such that, in $(-\varepsilon_0, \varepsilon_0)^g \times (-\varepsilon_1,\varepsilon_1)^{n-g}$,
$$F({\bf r}_0,{\bf r}_1)={\bf 0} \Leftrightarrow {\bf r}_0=g({\bf r}_1).$$

For each $l=0,\ldots,g$ such that $n_{l}> 1$ and $k=1,\ldots,n_{l}-1$, we want to construct a parametrized straight line 
$w_{j,k}:[0,T]\mapsto \overline{\Omega}$ such that $w_{j,k}(t)=w_{j,k}+r_{j,k}(t){\bf n}_{j,k}$, $0\le r_{j,k}(t) < \varepsilon_1$ and $u_j(w_{j,k}(t))=1-\frac{t}{n_j-1}$.
 Recall that
$u_{j}(w)=\delta_{jl}$ for $w\in \gamma_{l}$. Hence by the maximum principle of harmonic function, we have $u_j(z)<1$ for all $z \in \Omega$. It follows from Lemma \ref{lem:Hopf} that $\frac{\partial u_j}{\partial n}(w)\neq 0$ for any $w\in\gamma_{j}$. This implies that $u_j$ has no critical points on $\gamma_{j}$. Hence the level curve in $\Omega$ defined by $u_j(z)=1-\frac{t}{n_j-1}$ is a Jordan curve close to $\gamma_j$ for all sufficiently small positive $t$. We then define $w_{j,k}(t)=w_{j,k}+r_{j,k}(t){\bf n}_{j,k}$ to be the point on this level curve that intersects the normal line at $w_{j,k}$. We also need to define the parametrized straight line $w_{0,0}(t)=w_{0,0}+r_{0,0}(t){\bf n}_{0,0}$. 
Note that $\sum_{k=1}^{n_{j}-1} (1-u_{j}(w_{j,k}(t))=t$ tends to $0$ as $t\searrow 0$.  On the other hand, $u_j(w_{0,0})=0$ and $u_j(z)>0$ in $\Omega$. By Lemma \ref{lem:Hopf}, $\frac{\partial u_j}{\partial n}(w_{0,0})< 0$. Let $r_{0,0}(t)=\frac{1}{\alpha}t$ where $0<\alpha < |\frac{\partial u_j}{\partial n}(w_{0,0})|$ for each $1 \le j \le g$. Then the derivative of $u_j(w_{0,0}+r_{0,0}(t){\bf n}_{0,0})$ at $t=0$ is strictly greater than $1$. Since $u_{j}(w_{0,0}(0))=0=\sum_{k=1}^{n_{j}-1} (1-u_{j}(w_{j,k}(0)))$, we then have  
\begin{equation}u_{j}(w_{0,0}(t))>\sum_{k=1}^{n_{j}-1} (1-u_{j}(w_{j,k}(t))\label{eqn:lower bound}\end{equation}
for each $j=1,\ldots,g$ and $t\in(0,T]$ for sufficiently small positive $T$. 

Then for these $n-g$ parametrized lines $w_{0,0}(t)$ and $w_{l,k}(t)$, by the implicit function theorem, there  exist corresponding unique $w_{1,0}(t),\ldots ,w_{g,0}(t)$ such that 

\[\sum_{l=0}^{g}\sum_{k=0}^{n_{l}-1} u_{j}(w_{l,k}(t))=n_{j}.\]

We claim that all the $n$ points $w_{l,k}(t)$ are contained in $\Omega$. We only need to show that $w_{j,0}(t)\in\Omega$ for each $j=1,\ldots, g$. We first notice that by (\ref{eqn:lower bound}), we have

\begin{align*} n_{j}&=\sum_{l=0}^{g} \sum_{k=0}^{n_{l}-1} u_{j}(w_{l,k}(t))\\&\ge
\sum_{i=1}^{g} u_{j}(w_{i,0}(t))+u_{j}(w_{0,0}(t))+\sum_{k=1}^{n_{j}-1} u_{j}(w_{j,k}(t))
\\&=\sum_{i=1}^{g} u_{j}(w_{i,0}(t))+(n_j-1)-(n_j-1)+u_{j}(w_{0,0}(t))+\sum_{k=1}^{n_{j}-1} u_{j}(w_{j,k}(t))
\\&>\sum_{i=1}^{g} u_{j}(w_{i,0}(t))+(n_j-1).
\end{align*}
Hence, 
\begin{equation}
1>\sum_{i=1}^{g} u_{j}(w_{i,0}(t)) 
\label{eq:upper bound}\end{equation}
 for each $1\le j \le g$.

Suppose (for a contradiction) that for some $1\le i_1 \le g$, $w_{i_1,0}(t)\notin \Omega$ for all $t\in[0,T]$, then $u_{i_1}(w_{i_1,0}(t))\ge 1$. Since $\sum_{i=1}^{g} u_{j}(w_{i,0}(t)) <1$, we must have some different $1\le i_2 \le g$ such that $u_{i_1}(w_{i_2,0}(t))\le 0$ and hence $w_{i_2,0}(t) \notin \Omega$. In fact $w_{i_1,0}(t)$ and $w_{i_2,0}(t)$ are inside the disks bounded by the boundary circles $\gamma_{i_1}$ and $\gamma_{i_2}$ respectively. Let us assume that $i_1=1$ and $i_2=2$. If we define $u_A$ to be the harmonic measure of $\Omega$ with respect to $A \subset \partial\Omega$. Then $u_1+u_2=u_{\gamma_1\cup\gamma_2}$ on $\Omega$ because the two harmonic functions $u_1+u_2$ and $u_{\gamma_1\cup\gamma_2}$ have the same boundary values on $\partial\Omega$. Like what we have done before, we will also extend $u_{\gamma_1\cup\gamma_2}$ to the neighborhood $W$ of $\overline{\Omega}$ so that we still have $u_1+u_2=u_{\gamma_1\cup\gamma_2}$ on $W$. For this extended $u_{\gamma_1\cup\gamma_2}$, we know that $u_{\gamma_1\cup\gamma_2}(z)\ge 1$ for any $z \in W$ and inside the closed disks bounded by $\gamma_1$ and $\gamma_2$. Hence we have  $u_{\gamma_1\cup\gamma_2}(w_{1,0}(t))\ge 1$ and $u_{\gamma_1\cup\gamma_2}(w_{2,0}(t))\ge 1$ for $t\in(0,T]$. In general, we can also define harmonic measures for any union of the boundary components $\gamma_j$. Each such harmonic measure can then be extended across the boundary components so that the extension is greater than $1$ if it passes a boundary on which the harmonic measure is $1$
and less than zero if it passes a boundary on which the harmonic measure is $0$.

Summing up (\ref{eq:upper bound}) for $j=1$ and $2$, we have \begin{equation}
2>\sum_{i=1}^{g} u_{\gamma_1\cup\gamma_2}(w_{i,0}(t)).  \label{eq:upperbound2}\end{equation}
This implies that there must be some $i_3$ different from $1$ and $2$ such that $u_{\gamma_1\cup\gamma_2}(w_{i_3,0}(t))< 0$ and we may simply assume $i_3=3$. Then $w_{3,0}(t)$ will be inside the open disk bounded by $\gamma_3$. Continuing inductively, we conclude that for each $1\le j \le g$, $w_{j,0}(t)$ is inside the open disk bounded by $\gamma_j$ so that $u_0(w_{j,0}(t))<0$ for any $1\le j \le g$.

Now summing up (\ref{eq:upper bound}) from $j=1$ to $g$, we have
\begin{align*} g&>\sum_{i=1}^{g} \sum_{j=1}^{g} u_{j}(w_{i,0}(t))=\sum_{i=1}^{g} (1-u_{0}(w_{i,0}(t))=g-\sum_{i=1}^{g} u_{0}(w_{i,0}(t))>g\end{align*}
which is a contradiction. Here, we have used the fact that $u_0+\cdots+u_g\equiv 1$ on $\Omega$ (because the harmonic function $u_0+\cdots+u_g$ has boundary value $1$ on $\partial\Omega$). 

Therefore for any $t\in(0,T]$, the points $w_{l,k}(t)\in\Omega$ for $l=0,\ldots,g$ and $k=0,\ldots,n_{l}-1$ satisfy the condition (\ref{eq:condition1}). Lemma \ref{lem:allpoints} implies that there exists a proper holomorphic map $F_{t}\in\mathcal{P}(\Omega)$ with $F_{t}(p)=0$, $F_{t}(1)=1$ such that $F_{t}(w_{l,k}(t))$ does not depend on $l$ and $k$.
By the Schwarz reflection principle, each $F_t$ can be extended to a holomorphic map from a domain $W$ containing $\overline{\Omega}$ into some compact set $V \subset \mathbb{C}$. Then $\{F_{\frac{1}{n}}:n\in\mathbb{N}, \frac{1}{n}\in(0,T']\}$ forms a normal family on $W$ by Montel's Theorem. By passing to a convergent subsequence $f_{n}$, we can define $f:W \to V$ to be the local uniform limit of $f_{n}$ in $W$ as $n\rightarrow \infty$. Then $f$ is non-constant holomorphic in $\Omega$ and $f(p)=0$. Notice that for any two distinct $w_{l,k}$ and $w_{l',k'}$, we have $f(w_{l,k})=f(w_{l',k'})$. Here, we use the fact that $|f(w_{l,k})-f(w_{l',k'})|\le |f(w_{l,k})-f_n(w_{l,k}(\frac{1}{n}))|+|f_n(w_{l',k'}(\frac{1}{n}))-f(w_{l',k'})|$ and $|f(w_{l,k})-f_n(w_{l,k}(\frac{1}{n}))| \to 0$ as $n \to \infty$ because $f_n \to f$ uniformly in any compact subset of $W$ containing $w_{l,k}$. Since  $|f_n(w_{l,k})|=1$, we deduce that $f( w_{l,k})$ takes the same value in $\partial \mathbb{D}$ for any $w_{l,k}$. By composing with a rotation, we can assume that $f( w_{l,k} )=1$ for any $l$ and $k$. Then, the set of zeros of $f$ are the limits of the zeros of $f_{n}$ as $n\rightarrow \infty$. By the continuity of harmonic measure, the zeros of $f$ also satisfy the condition (\ref{eq:condition1}) and, by continuity of $\eta(w,\cdot)$, we deduce that $f$ takes the form in Theorem \ref{thm:expression} and hence $f$ is proper.
\\

We now prove the uniqueness part. Suppose that there exists two proper holomorphic maps $f_{1}, f_{2}\in\mathcal{P}_{\nu}(\Omega)$ that satisfy conditions (i), (ii) and (iii) from the statement in the theorem, i.e. $f_{1}(p)=f_{2}(p)=0$; $f_{1}(w_{l,k})=1=f_{2}(w_{l,k})$ for all $w_{l,k}$ and for each $l=0,\ldots,g$ such that $n_{l}> 1$ and $k=1,\ldots,n_{l}-1$, 
\begin{equation}
\frac{|f_1'(w_{0,0})|}{|f_1'(w_{l,k})|}=\lambda_{l,k}=\frac{|f_2'(w_{0,0})|}{|f_2'(w_{l,k})|}.\label{eq:derivative}
\end{equation}

Since $f_i'(w_{l,k})\neq 0$ for all $w_{l,k,}$, $f_i$ maps  some neighborhood $U_{l,k}$ of $w_{l,k,}$ conformally onto a fixed neighborhood $U$ of $1$. Let $w_{0,0}^i(t)=w_{0,0}+t {\bf n_{0,0}}$ for $t\in [0,T]$ be a straight line segment $L$ in $U_{0,0}$. As $f_i$ is angle preserving at $w_{0,0}$, $f^i(L)$ is orthogonal to $\partial\mathbb{D}$ at $1$. Then $f_i^{-1}(f_i(L))$ consists of disjoint curves $w_{l,k}^i(t)$, each  with an endpoint at some $w_{l,k}$ when $t=0$ and is orthogonal to the tangent line of $\gamma_l$ at $w_{l,k}$. It follows that 
$f_i(w_{l,k}(t))=f_i(w_{0,0}(t))$ and hence 
$|f_i(w_{l,k})w_{l,k}'(0)|=|f_i(w_{0,0})w_{0,0}'(0)|=|f_i(w_{0,0})|$. If we take $T$ to be sufficiently close to $0$, then each $w_{l,k}^i(t)$ is a line segment of the form $w_{l,k}^i(t)=w_{l,k}^i+r_{l,k}^i(t){\bf n_{l,k}}$. From (\ref{eq:derivative}), for $i=1,2$, we have  $w_{l,k}^i(t)=w_{l,k}+\lambda_{l,k}t {\bf n_{l,k}}$ for each $l=0,\ldots,g$ such that $n_{l}> 1$ and $k=1,\ldots,n_{l}-1$. Since $f_i(w_{l,k}^i(t))=f_i(w_{0,0}(t))$ for all $l$ and $k$, by the first part of Lemma \ref{lem:allpoints}, we have \[\sum_{l=0}^{g}\sum_{k=0}^{n_{l}-1} u_{j}(w_{l,k}^i(t))-n_{j}=0\]
where $w_{0,0}^1(t)=w_{0,0}+t {\bf n_{0,0}}=w_{0,0}^2(t)$ and $w_{l,k}^1(t)=w_{l,k}+\lambda_{l,k}t {\bf n_{l,k}}=w_{l,k}^2(t)$ for each $l=0,\ldots,g$ such that $n_{l}> 1$ and $k=1,\ldots,n_{l}-1$. Like what we did previously, we apply the implicit function theorem to (\ref{eqn:implicitcondition}) and conclude from the uniqueness part of the theorem that $w_{j,0}^1(t)=w_{j,0}^2(t)$ for $j=1,\ldots g$. It follows that $v_{l,k}:=w_{l,k}^1(T)=w_{l,k}^2(T)$ for any $l$ and $k$ and hence for some $\zeta_{1},\zeta_{2}\in\mathbb{D}$, $f_{1}(v_{l,k})=\zeta_{1}$ and $f_{2}(v_{l,k})=\zeta_{2}$ for all $l$ and $k$. 

Now, for $i=1,2$, let $\mu_{i}$ be a M\"{o}bius transformation that fixes $\mathbb{D}$ and maps $\zeta_{i}$ to $0$. So $\mu_{1}\circ f_{1}$ and $\mu_{2}\circ f_{2}$ are proper holomorphic maps in $\mathcal{P}_{n}(\Omega) $ and  have the same zeros. By Theorem \ref{thm:expression}, we have $\mu_{1}\circ f_{1}=\lambda\mu_{2}\circ f_{2}$ for some $|\lambda|=1$. But $f_{1}(p)=f_{2}(p)=0$ and $f_{1}(w_{0,0})=f_{2}(w_{0,0})=1$. Thus $\mu_{1}(0)=\lambda\mu_{2}(0)$ and $\mu_{1}(1)=\lambda\mu_{2}(1)$. This implies that $\mu_{1}=\lambda\mu_{2}$ and so $f_{1}=f_{2}$. 
 \end{proof}
We now prove Corollary \ref{cor:parameter}. 
\begin{proof}[Proof of Corollary \ref{cor:parameter}]
In this proof, we will prove the result for the map from $\Omega  \times \mathbb{S}^{n} \times \mathbb{R}_{+}^{n-g-1}$ onto  $\mathcal{P}_{\nu}(\Omega)$. The other case follows in a similar way.

Theorem \ref{thm:parameter} proves the existence of the surjective map $H: \Omega  \times \mathbb{S}^{n} \times \mathbb{R}_{+}^{n-g-1} \to \mathcal{P}_{\nu}(\Omega)$. To see that $H$ is $(n\cdot n_0 !\cdots n_g !)$-to-$1$, note that we can apply a permutation to the specified points on each boundary to get the same proper map under Theorem \ref{thm:parameter}. Since there are $n_{l}$ points specified on the boundary component $\gamma_{l}$ for each $l=0,\ldots g$, there are in total $(n_0 !\cdots n_g !)$ such permutations. The uniqueness part of Theorem \ref{thm:parameter} guarantees that, for a fixed $p\in\Omega$, these $(n_0 !\cdots n_g !)$ parameters are the only parameters that give the same proper map with a zero at $p$. In addition, any one of the $n$ zeros of the proper map could be given as the $p$ in the statement of Theorem \ref{thm:parameter}. Hence, we deduce that $H$ is $(n\cdot n_0 !\cdots n_g !)$-to-1.

To obtain the continuity of $H$, take any ${\bf w}\in \Omega  \times \mathbb{S}^{n} \times \mathbb{R}_{+}^{n-g-1}$ and a sequence $\{{\bf w}_{k}\} \subset \Omega  \times \mathbb{S}^{n} \times \mathbb{R}_{+}^{n-g-1}$ such that ${\bf w}_{k}\rightarrow {\bf w}$ as $k\rightarrow\infty$. For $k=1,2,\ldots$, let 
\[{\bf w}_{k}=(p^{k},w^{k}_{0,0},\ldots,w^{k}_{g,n_{g}-1},\lambda^{k}_{0,1},\ldots,\lambda^{k}_{0,n_{0}-1},\ldots,\lambda^{k}_{g,1},\ldots, \lambda^{k}_{g,n_{g}-1})\] and 
\[{\bf w}=(p,w_{0,0},\ldots,w_{g,n_{g}-1},\lambda_{0,1},\ldots,\lambda_{0,n_{0}-1},\ldots,\lambda_{g,1},\ldots, \lambda_{g,n_{g}-1}).\] 
Now let $f_{k}$ be the proper holomorphic map in $\mathcal{P}_{\nu}(\Omega)$ which satisfies Theorem \ref{thm:parameter} for the point ${\bf w}_{k}$ and let $f$ denote the proper holomorphic map in $\mathcal{P}_{\nu}(\Omega)$ which satisfies Theorem \ref{thm:parameter} for the point ${\bf w}$. As we have mentioned in Section 1, we can extend $f_k$ and $f$ to holomorphic maps in a neighborhood $U$ of $\overline{\Omega}$ so that they are still bounded on $U$. Since $\{f_k\}$ is a family of bounded holomorphic maps on $U$, by Montel's Theorem, there exists a subsequence $\{f_{k_i}\}$ such that $f_{k_i}$ converges to a holomorphic map $h$ locally uniformly on $U$.   Since ${\bf w}_{k_i} \to {\bf w}$, and $f_{k_i} \to h$ locally uniformly, by Hurwitz's theorem, we must have $h(p)=0$ and $h^{-1}(1)=\{w_{0,0},\ldots,w_{g,n_{g}-1}\}$. In addition, by Weierstrass theorem, we have $f_{k_i}' \to h'$ locally uniformly on $U$ and together with the fact that each $f_{k_i}$ satisfies condition (iii) in Theorem \ref{thm:parameter} for the set $\{w^{k_i}_{0,0},\ldots,w^{k_i}_{g,n_{g}-1}\}$, we have \[\frac{|h'(w_{0,0})|}{|h'(w_{l,k})|}=\lambda_{l,k}.\]
By the uniqueness part of Theorem \ref{thm:parameter}, we must have $f=h$ and therefore $f_{k_i} \to f$ locally uniformly on $\Omega$. 

Suppose that $H$ is not continuous at some ${\bf w} \in  \Omega \times \mathbb{S}^{n} \times \mathbb{R}_{+}^{n-g-1}$. Since $\mathcal{P}_{\nu}(\Omega)$ is equipped with the  topology of local uniform convergence, we can find a sequence
$\{{\bf w}_{k}\} \subset \Omega  \times \mathbb{S}^{n} \times \mathbb{R}_{+}^{n-g-1}$ such that ${\bf w}_{k}\rightarrow {\bf w}$ and some compact $K \subset \Omega$, such that for the sup norm $||\cdot||_K$,
\[||H({\bf w}_{k})-H({\bf w})||_K =||f_k-f||_K > \varepsilon >0\]
for all $k \in \mathbb{N}$. On the other hand, we have shown above that we can find a subsequence $f_{k_i}$ which converges to $f$ locally uniformly on $\Omega$ which is a contradiction.   
\end{proof}
\section{Finite Blaschke products and proper holomorphic maps}
In this section, we will establish a link between two function spaces  $\mathcal{P}(\Omega)$ and $\mathcal{P}(\mathbb{D})$ where $\Omega$ is a circular domain. Recall that by a famous result of Fatou \cite{fatou}, every proper holomorphic map in $\mathcal{P}(\mathbb{D})$ is a finite Blaschke product. It is a standard result that a finite Blashcke product is uniquely determined by its zeros and image of 1. Also, finite Blaschke products satisfy the following properties: for any $f,g\in\mathcal{P}(\mathbb{D})$
\begin{enumerate}
\item[(i)] $f\cdot g\in\mathcal{P}(\mathbb{D})$ with $\deg (f\cdot g) = \deg f +\deg g$;
\item [(ii)]$ f\circ g\in\mathcal{P}(\mathbb{D})$ with $\deg (f\circ g) = \deg f \cdot \deg g$.
\end{enumerate}
A similar result holds for $\mathcal{P}(\Omega)$.
\begin{proposition}\label{prop:embed}
For any $f,g\in\mathcal{P}(\Omega)$ and $h\in\mathcal{P}(\mathbb{D})$,
\begin{enumerate}
\item[(i)] $f\cdot g\in\mathcal{P}(\Omega)$ with $\deg (f\cdot g) = \deg f + \deg g$ and 
$\left. \deg (f \cdot g \right|_{\gamma_{j}}) = \left. \deg (f\right|_{\gamma_{j}}) + \left.\deg (g\right|_{\gamma_{j}})$ for each $j=0,\ldots,g$;
\item [(ii)]$ h\circ f\in\mathcal{P}(\Omega)$ with $\deg (h\circ f) = \deg h \cdot \deg f$ and 
$\deg \left(  \left. h \circ f\right|_{\gamma_{j}} \right) =  \deg h \cdot \left(\deg \left. f\right|_{\gamma_{j}} \right)$ for each $j=0,\ldots,g$. 
\end{enumerate}

\end{proposition}
\begin{proof}
Suppose that $f,g \in \mathcal{P}(\Omega)$ such that $f$ has degree $n$, boundary degree $(n_{0},n_{1},\ldots,n_{g})$, and zeros $p_{1},\ldots,p_{n}\in\Omega$; $g$ has degree $m$, boundary degree $(m_{0},m_{1},\ldots,m_{g})$, and zeros $q_{1},\ldots,q_{m}\in\Omega$.
Hence by condition (\ref{eq:condition1}),
\[     \sum_{k=1}^n u_j (p_k) = n_j \mbox{ and } \sum_{k=1}^m u_j (q_k) = m_j  \mbox{ for all $j=1,\dots , g$}.\]
To prove (i), we note that
 $p_{1},\ldots,p_{n},q_{1},\ldots, q_{m}$ satistfy 
\[\left(\sum_{k=1}^n u_j (p_k)\right)+ \left(\sum_{k=1}^m u_j (q_k)\right) =n_{j}+m_{j} 
     \qquad
     \mbox{for all $j=1,\dots , g$}.\]
Condition (\ref{eq:condition1}) implies that there exists a proper holomorphic map, $F\in\mathcal{P}(\Omega)$, with degree $n+m$ and boundary degree $(n_{0}+m_{0},\ldots, n_{g}+m_{g})$ with these zeros. As in the proof of the uniqueness part of Theorem \ref{thm:expression}, each of $F$, $f$ and $g$ extend to holomorphic maps from the compact Riemann surface $\overline{\Omega^{*}}/\Theta$ onto the Riemann sphere. Then $G=\dfrac{F}{f\cdot g}$ also defines a holomorphic map on the compact Riemann surface $\overline{\Omega^{*}}/\Theta$ with no zeros or poles. Hence $G$ is constant.
Thus $f\cdot g= \lambda F$ for some $\lambda$ with $|\lambda|=1$ and so $f\cdot g\in\mathcal{P}(\Omega)$.

To prove (ii), we suppose $h \in\mathcal{P}(\mathbb{D})$ is of degree $N$ and has zeros $a_{0},\ldots, a_{N}$. 
For $l=1,\ldots,N$, let $p_{1,l},\ldots,p_{n,l}$ denote the preimages of $a_{l}$ under $f$. By Lemma \ref{lem:allpoints}, for all $l=1,\ldots,N$,
\[     \sum_{k=1}^n u_j (p_{l,k}) = n_j
     \qquad
     \mbox{for all $j=1,\dots , g$}.\]
Hence,
\[     \sum_{l=1}^{N}\sum_{k=1}^n u_j (p_{l,k}) = N\cdot n_j
     \qquad
     \mbox{for all $j=1,\dots , g$}.\]
Condition (\ref{eq:condition1}) implies that there exists a proper holomorphic map, $\widetilde{F}\in\mathcal{P}(\Omega)$ with degree $Nn$, boundary degree $(Nn_{0},\ldots, Nn_{g})$ with zeros at $p_{l,k}$ for each $l=1,\ldots,N$ and $k=1,\ldots,n$. Proceeding with the same argument as above, we deduce that $H=\dfrac{\widetilde{F}} {h\circ f}$ is constant. This implies that $h\circ f=\lambda'\widetilde{F}$ for some $\lambda'$ with $|\lambda'|=1$. Thus $h\circ f\in\mathcal{P}(\Omega)$.
\end{proof}
Let $\Phi:\mathcal{P}(\Omega)\mapsto \mathcal{P}(\mathbb{D})$ send a proper holomorphic map in $\mathcal{P}(\Omega)$ to a finite Blaschke product with the same zeros and the same image of $1$. It is easy to see that this map is well-defined and injective. However, in the case where the infinite product expression  (\ref{eq:def_skpf}) holds, the map $\Phi:\mathcal{P}(\Omega)\mapsto \mathcal{P}(\mathbb{D})$ can be made explicit and moreover preserves the multiplicative semigroup structure given in Proposition \ref{prop:embed}. This is given as the following result.
\begin{theorem}
\label{thm:embedding}
For any $f\in\mathcal{P}(\Omega)$, let $B\in\mathcal{P}(\mathbb{D})$ be a finite Blaschke product with the same zeros as $f$ such that $B(1)=1$. Assume that  (\ref{eq:def_skpf}) holds, then we can write $f$ as 
\begin{equation}f(z)=\lambda \exp \left( -2 \pi i \sum\limits_{j=1}^{g} n_j v_j (z) \right)\prod_{\theta\in\Theta} \dfrac{B(\theta(z))}{B(\theta(1))}
\label{eqn:embed}\end{equation}
for some constant $\lambda$ with $|\lambda|=1$.
If we denote this $f$ by $f_{B}$, and let $\Phi$ map $f_{B}\in\mathcal{P}(\Omega)$ to $B \in \mathcal{P}(\mathbb{D})$, then 
for any $f, g \in \mathcal{P}(\Omega)$
\[\Phi(f\cdot g)(z)=\Phi(f)\cdot\Phi(g).\]
\end{theorem}
\begin{proof}
The expression (\ref{eqn:embed}) follows directly by combining Theorem \ref{thm:expression} and Proposition \ref{lem:embedding}. 

Now suppose that $f_{B_{1}}, f_{B_{2}}\in\mathcal{P}(\Omega)$ for some finite Blaschke products $B_{1}$ and $B_{2}$. Proposition \ref{prop:embed}(i) implies that $f_{B_{1}}\cdot f_{B_{2}}\in\mathcal{P}(\Omega)$ and hence we can write $f_{B_{1}}\cdot f_{B_{2}} = f_{B}$ for some finite Blaschke product $B$.
From a comparison of the zeros and the image of $1$, we deduce that $B= B_{1}\cdot B_{2}$. This establishes the desired result.
\end{proof}

\bibliographystyle{amsplain}
\bibliography{reference}
\label{pg:reference}

\end{document}